\documentclass[a4paper,11pt]{article}
\usepackage{latexsym}
\usepackage{amsmath}
\usepackage{amsthm}

\usepackage{graphicx}
\usepackage{txfonts}
\usepackage{bm}
\usepackage{color}	
\usepackage{epic,eepic}

\newcommand{\BLU}[1]{{\color{blue}#1}}

\renewcommand{\BLU}[1]{{\color{black}#1}}

\usepackage{geometry}
\numberwithin{equation}{section}

\newtheorem{theorem}{Theorem}[section]
\newtheorem{proposition}[theorem]{Proposition}
\newtheorem{corollary}[theorem]{Corollary}
\newtheorem{lemma}[theorem]{Lemma}
\newtheorem{claim}[theorem]{Claim}

\newtheorem{remark}[theorem]{Remark}

\newcommand{\OMIT}[1]{{\bf [OMIT:} #1 \ {\bf --- end OMIT] }}  
   \renewcommand{\OMIT}[1]{}            

\newcommand{\ZZ}{{\bf Z}}

\newcommand{\odotZ}{\overset{..{.}.}}

\newcommand{\Proof}{\noindent {\bf Proof.  }}

\newcommand{\nat}{\sp{\natural}}
\newcommand{\finbox}{\hspace*{\fill}$\rule{0.17cm}{0.17cm}$}
\newcommand{\FBOX}{\hspace*{\fill}$\rule{0.17cm}{0.17cm}$}
\newcommand{\BB}{\hspace*{\fill}$\rule{0.17cm}{0.17cm}\ \rule{0.17cm}{0.17cm}$}

\usepackage{lineno}
\newcommand*\patchAmsMathEnvironmentForLineno[1]{
  \expandafter\let\csname old#1\expandafter\endcsname\csname #1\endcsname
  \expandafter\let\csname oldend#1\expandafter\endcsname\csname end#1\endcsname
  \renewenvironment{#1}
     {\linenomath\csname old#1\endcsname}
     {\csname oldend#1\endcsname\endlinenomath}}
\newcommand*\patchBothAmsMathEnvironmentsForLineno[1]{
  \patchAmsMathEnvironmentForLineno{#1}
  \patchAmsMathEnvironmentForLineno{#1*}}
\AtBeginDocument{
\patchBothAmsMathEnvironmentsForLineno{equation}
\patchBothAmsMathEnvironmentsForLineno{align}
\patchBothAmsMathEnvironmentsForLineno{flalign}
\patchBothAmsMathEnvironmentsForLineno{alignat}
\patchBothAmsMathEnvironmentsForLineno{gather}
\patchBothAmsMathEnvironmentsForLineno{multline}
}

\begin{document}

\title{
Fair Integral Submodular Flows
}

\author{Andr\'as Frank%
\thanks{MTA-ELTE Egerv\'ary Research Group,
Department of Operations Research, E\"otv\"os University, P\'azm\'any
P.~s.~1/c, Budapest, Hungary, H-1117. 
e-mail:  {\tt frank@cs.elte.hu}. 
ORCID: 0000-0001-6161-4848.
The research was partially supported by the Hungarian
Scientific Research Fund - OTKA, No. NKFIH-FK128673.
}
 \ \ and \ 
{Kazuo Murota%
\thanks{
\BLU{
Faculty 
}%
of Economics and Business Administration,
Tokyo Metropolitan University, Tokyo 192-0397, Japan;
\BLU{
Currently at
The Institute of Statistical Mathematics,
Tokyo 190-8562, Japan;
}%
e-mail:  {\tt murota@tmu.ac.jp}.
ORCID: 0000-0003-1518-9152.
The research was supported by 
JSPS KAKENHI Grant Number JP20K11697. 
 }}}


\date{December 2020 / June 2022}

\maketitle

\begin{abstract}
Integer-valued elements of an integral submodular
flow polyhedron $Q$ are investigated which are decreasingly minimal
(dec-min) in the sense that their largest component is as small as
possible, within this, the second largest component is as small as
possible, and so on.  As a main result, we prove that the set of
dec-min integral elements of $Q$ is the set of integral elements of
another integral submodular flow polyhedron arising from $Q$ by
intersecting a face of $Q$ with a box.  Based on this description, we
develop a strongly polynomial algorithm for computing not only a
dec-min integer-valued submodular flow but even a cheapest one with
respect to a linear cost-function.  A special case is the problem of
finding a strongly connected (or $k$-edge-connected) orientation of a
mixed graph whose in-degree vector is decreasingly minimal.
\end{abstract}

\\ \qquad
Polynomial algorithm.







\medskip



\section{Introduction} \label{Intro}

For an integral polyhedron $Q$, let $\odotZ{Q}$ denote the set of
integral elements of $Q$.  
We are interested in finding the most fair element $z$ of $\odotZ{Q}$,
where fairness wants to reflect the intuitive feeling that 
the components of $z$ are distributed as equitably as possible.  
There may be several ways to formally capture fairness.  
For example, if the difference of the largest and the
smallest components of $z$ is minimum, then $z$ is felt rather fair.
The square sum of the components is another, more global measure for fairness.  
In this work we are interested in a third natural possibility:  
decreasing minimality.  Let us call an element 
$z\in \odotZ{Q}$ decreasingly minimal (dec-min) if the largest component of $z$ 
is as small as possible, within this, the second largest component
(possibly with the same value as the largest component) is as small as possible, 
within this, the third largest component of $z$ is small as possible, 
and so on.  
Actually, it is more convenient to work with the following slightly more general concept.  
Let $F$ be a specified subset of the coordinates.  
We say that $z\in \odotZ{Q}$ is {\bf decreasingly minimal on} $F$ 
(or {\bf $F$-dec-min} for short) 
if the restriction of $z$ to $F$ is decreasingly minimal.  
For the trivial special case $F=\emptyset $, the set of $F$-decmin elements of 
$\odotZ{Q}$ is $\odotZ{Q}$ itself.


In \cite{Frank-Murota.A} and \cite{Frank-Murota.B}, we characterized
dec-min elements of an M-convex set (which is, by an equivalent
definition, the set $\odotZ{B}$ of integral elements of an integral base-polyhedron $B$ 
\cite{Murota98a,Murota03}) 
and described 
a strongly polynomial algorithm to compute a cheapest dec-min element
with respect to a linear cost function.  
In \cite{Frank-Murota.C}, an analogous investigation was carried out for network flows.  
In the present work, we consider decreasingly minimal integer-valued
submodular flows, a common generalization of network flows and base-polyhedra.

Submodular flows were introduced by Edmonds and Giles \cite{Edmonds-Giles} 
while the term itself was suggested by Zimmermann \cite{Zimmermann82b}.  
The notion became a standard tool in discrete optimization, 
see the books \cite{Frank-book,Fujishige-book,Murota03,Schrijverbook}.  
Originally, Edmonds and Giles used submodular functions for the definition but
supermodular functions could equally well be used, and in applications
one often needs supermodular functions.  
Therefore, in the present paper, 
we replace the term submodular flow by {\bf base-flow} 
(but emphasize that these are the same). From this point of view, the term
base-polyhedron is adequate since it refers to neither submodular nor supermodular functions.  
Our newly suggested term base-flow is intended to capture this neutrality.  
In Section~\ref{basic-prop}, we
recall the formal definitions of base-polyhedra and base-flows, along
with their basic properties.

\subsection{Main results} 
\label{main-results}

One of our main goals is to provide a complete description of the set
of $F$-dec-min elements of $\odotZ{Q}$ for an integral base-flow polyhedron $Q$.  
Roughly, the theorem states that this set is the set of integral elements 
of a base-flow polyhedron obtained from $Q$ by
intersecting a face of $Q$ with a box which is `narrow' on $F$.
In what follows, ${\bf Z}$ denotes the set of integers while 
$\overline{\bf Z}:= {\bf Z} \cup \{+\infty \}$ and 
$\underline{\bf Z}:= {\bf Z} \cup \{-\infty \}$.  
Throughout we use ${\bf R}$ to denote the set of reals.


\begin{theorem}
\label{MAIN} 
Let $D=(V,A)$ be a digraph endowed with integer-valued 
lower and upper bound functions $f:A\rightarrow \underline{\bf Z}$ 
and $g:A\rightarrow \overline{\bf Z}$ 
for which $f\leq g$.  
Let $B$ be an integral base-polyhedron for which the base-flow polyhedron
$Q=Q(f,g;B)$ is non-empty.  Let $F\subseteq A$ be a specified subset
of edges such that both $f$ and $g$ are finite-valued on $F$.  
Then there exists a face $B\sp{\triangledown}$ of $B$ and there exists a pair
$(f\sp{*},g\sp{*})$ of integer-valued bounding functions on $A$ with
$f\leq f\sp{*} \leq g\sp{*} \leq g$ 
such that an element $z\in \odotZ{Q}$ is 
$F$-dec-min if and only if $z\in \odotZ{Q}(f\sp{*},g\sp{*};B\sp{\triangledown})$.  
Moreover, $0\leq g\sp{*}(e)-f\sp{*}(e)\leq 1$ for every $e\in F$.  
\end{theorem}

The proof will be prepared in Sections \ref{flow3} and \ref{minmaxL},
and completed in Section~\ref{decmin}.  
It should be noted that in the general case, 
when the finiteness of $f$ and $g$ on $F$ is not assumed, 
it may be the case that no $F$-dec-min element exists at all.
In Section~\ref{vandecmin} we describe a characterization for the
existence of $F$-dec-min integral feasible base-flow, and show that
Theorem~\ref{MAIN} extends to this case, as well.

It is a known and easy property \cite{FrankP6} that the intersection
of two (integral) $g$-polymatroids is an (integral) base-flow polyhedron.  
Therefore Theorem~\ref{MAIN} can be specialized to the following.


\begin{corollary} 
\label{M2} 
Let $Q_{1}$ and $Q_{2}$ be two 
bounded integral g-polymatroids in ${\bf R} \sp{S}$ 
for which their intersection $Q$ is non-empty.  
Let $F\subseteq S$ be a specified subset of ground-set $S$.  
Then there exists a face 
$Q\sp{\triangledown}$ of $Q$ and there exists
an integral box $T=T(f\sp{*},g\sp{*})$ in ${\bf R} \sp{S}$ 
such that an element $z\in \odotZ{Q}$ is $F$-dec-min 
if and only if 
$z\in {\odotZ{Q\sp{\triangledown}}}\cap T$.  
Moreover, $0\leq g\sp{*}(s)-f\sp{*}(s)\leq 1$
for every element $s\in F$.
\end{corollary}


Recall that an M-convex set is the set of integral elements of 
an integral base-polyhedron.  
In Discrete convex analysis \cite{Murota98a,Murota03}, 
the set of integral elements of an integral
g-polymatroid is called an {\rm M}$\nat$-convex set
 (pronounce M-nat-convex or M-natural-convex).  
As a base-polyhedron is a special g-polymatroid, 
an M-convex set is {\rm M}$\nat$-convex.  
Furthermore, the intersection of two M-convex 
(resp., {\rm M}$\nat$-convex) sets is called an {\rm M}$_{2}$-convex 
(resp., {\rm M}$_{2}\nat$-convex) set.  
It was proved by Frank \cite{FrankP6} that 
the intersection of two integral g-polymatroids is an integral polyhedron, 
implying that an {\rm M}$_{2}\nat$-convex set is 
the set of integral elements of the intersection of two integral g-polymatroids.  
Therefore, Corollary~\ref{M2} can be interpreted as a
characterization of the set of $F$-dec-min elements of an {\rm M}$_{2}\nat$-convex set.

Our second main goal is to develop a strongly polynomial algorithm for
computing the bounding functions $f\sp{*}$ and $g\sp{*}$ in the theorem,
as well as the face $B\sp{\triangledown}$ of $B$.  
Once these data are available, with the help of a standard base-flow algorithm
\cite{FrankJ8,Frank-book}, one can compute an $F$-dec-min element of $\odotZ{Q}$.  
Even more, with the help of a minimum cost base-flow algorithm \cite{FrankJ15,FRZ89},
 a minimum cost $F$-dec-min element of
$\odotZ{Q}$ can also be computed in polynomial time, with respect to 
a linear cost-function $c:A\rightarrow {\bf R}$.  
Indeed, by Theorem~\ref{MAIN},  this latter problem is nothing but a minimum cost 
$(f\sp{*},g\sp{*})$-bounded base-flow problem, which can be solved in polynomial time. 
 (See, e.g.~\cite{FrankJ15,FRZ89}, and the book of
Schrijver \cite{Schrijverbook} (p.~1019)).

This approach gives rise to an algorithm in the special base-flow
problem when the goal is to find a (minimum cost) $k$-edge-connected
and in-degree-bounded orientation of a mixed graph for which the
in-degree vector is decreasingly minimal.


\begin{remark} \rm \label{RMconvmin}
The decreasing minimization problem 
is often related to minimization of a convex cost function.
For example, if $Q$ is a base-polyhedron,
an element of $Q$ is dec-min in $Q$
if and only if it is a square-sum minimizer of $Q$ \cite{Fujishige80,Fujishige-book}.
The corresponding statement is also true in its discrete version 
where $Q$ is an M-convex set \cite{Frank-Murota.A}.
Such equivalence between dec-minimality and 
square-sum minimality fails for network flows,
which is demonstrated by an example in \cite[Section~11.2]{Frank-Murota.C}.
It is also noted that
the dec-min problem for $Q \subseteq \ZZ\sp{n}$ (in general)
can be formulated as a separable convex function minimization,
as discussed in \cite[Section~3]{FM19partII}.
In our integral $F$-dec-min base-flow problem, we can take, for example,
a real-valued cost function
$\sum_{e \in F} |F|^{x(e)}$
for an integral base-flow $x$.
Here the function 
$\varphi(k) = |F|^{k}$, defined for all integers $k$,
is strictly convex in the sense of
$\varphi(k-1) + \varphi(k+1) > 2 \varphi(k)$ $(k \in \ZZ)$
and `rapidly increasing' in the sense of 
$\varphi(k+1) \geq  |F|  \  \varphi(k)$ $(k \in \ZZ)$.
Such convex formulation does not readily provide a strongly polynomial algorithm,
although it gives us structural insight.
$\bullet $ 
\end{remark}

\subsection{Basic notions and notation} \label{basic-prop}

Let $S$ be a finite ground-set.  
Two subsets $X$ and $Y$ of $S$ are {\bf intersecting} if $X\cap Y\not =\emptyset$. 
They are {\bf properly intersecting} if none of $X\cap Y$, $X-Y$, $Y-X$ is empty.
If, in addition, $S-(X\cup Y)$ is also non-empty, 
we speak of a {\bf crossing} pair.

For a vector $x\in {\bf R}\sp{S}$ or a function 
$x:S\rightarrow {\bf R}$, 
the set-function 
$\widetilde x:2\sp{S}\rightarrow {\bf R}$ 
is defined by 
\[
\widetilde x(Z):  = \sum [x(s):s\in Z] \quad (Z\subseteq S). 
\]
\noindent 
Such a function is modular in the sense that the modular equality 
\begin{equation} 
\widetilde x(X) + \widetilde x (Y)
 = \widetilde x(X\cap Y)+ \widetilde x (X\cup Y) 
\label{(modular)} 
\end{equation} 
\noindent 
holds for every pair $\{ X,Y \}$ of subsets of $S$.

Unless stated otherwise, we assume throughout for a set-function that
its value on the empty set is zero.  
Let $p:S\rightarrow \underline{\bf Z}$
be a set-function on $S$.  
When $p(S)$ is finite, the set-function $\overline{p}$ defined by 
\[
\overline{p}(X):  = p(S)- p(S-X)
\] 
\noindent 
is called the {\bf complementary set-function} or just the {\bf complement} of $p$.
Clearly, $\overline{p}(\emptyset )=0$ and $\overline{p}(S) = p(S)$, and the
complement of $\overline{p}$ is $p$ itself.

The supermodular inequality for $X,Y\subseteq S$ is as follows:  
\begin{equation}
 p(X) + p (Y) \leq p(X\cap Y)+ p (X\cup Y).  
\label{(supermod)} 
\end{equation}
\noindent 
The set-function $p$ is called {\bf fully supermodular} 
or just {\bf supermodular} if \eqref{(supermod)} 
holds for every pair of subsets $X,Y\subseteq S$.  
We say that $p$ is {\bf intersecting} 
({\bf crossing}) {\bf supermodular} if \eqref{(supermod)} holds 
for every intersecting (crossing) pair of subsets of $S$.

For a set-function $b:S\rightarrow \overline{\bf Z}$, 
the submodular inequality is as follows:  
\begin{equation} 
 b(X) + b (Y) \geq b(X\cap Y)+ b (X\cup Y). 
\label{(submod)} 
\end{equation}
\noindent 
Function $b$ is called fully (intersecting, crossing) submodular 
if \eqref{(submod)} holds for every (intersecting, crossing) pair of subsets.  
Clearly, the complement of a fully (crossing) supermodular function is fully (crossing) submodular.

With a set function $h$ with finite $h(S)$, we associate two polyhedra:  
\begin{align*}
B'(h)&:=\{x\in {\bf R} \sp{S}:  \widetilde x(S) = h(S) \
\hbox{and}\ \ \widetilde x(Z)\geq h(Z) \ \hbox{for all}\ \ Z\subset S
\}, 
\\
B(h) &:=\{x\in {\bf R} \sp{S}:  \widetilde x(S) = h(S) \
\hbox{and}\ \ \widetilde x(Z)\leq h(Z) \ \hbox{for all}\ \ Z\subset S
\}.  
\end{align*}
\noindent 
Obviously, $B(h)=B'(\overline{h})$ and $B'(h)=B(\overline{h})$.

For a (fully) supermodular function $p$ with finite $p(S)$, 
the polyhedron $B:=B'(p)$ is called a {\bf base-polyhedron}.  
We say that $B$ is a {\bf $0$-base-polyhedron} if $p(S)=0$, or equivalently,
$\widetilde x(S)=0$ for each $x\in B$.  
A base-polyhedron can also be described 
with the help of a submodular function $b$, namely, $B=B(b)$,
where $b$ and $p$ are complementary set-functions.  
The empty set, by convention, is also considered a base-polyhedron 
although it cannot be defined by a fully sub- (or supermodular) function 
(which is equivalent to saying that $B'(p)$ is never empty).

A basic property of non-empty base-polyhedra is that they uniquely
determine their bounding fully supermodular (or submodular) function.
Namely, 
\begin{equation} 
p(Z)=\min \{\widetilde x(Z):  x\in B\} 
\qquad (\hbox{or} \ \  
b(Z)=\max\{\widetilde x(Z):  x\in B\}).
\label{(pbdef)} 
\end{equation} 
\noindent 
This formula is particularly important from an algorithmic point of view
since, if an algorithm is developed 
for the case when the base-polyhedron $B$
is defined by a fully supermodular $p$, the same algorithm can be used
for an arbitrarily given base-polyhedron 
(defined, for example, by a crossing supermodular function),
provided that a subroutine is available to compute $p(Z)$ in \eqref{(pbdef)} 
for any input subset $Z\subseteq S$.

A rich overview of properties of base-polyhedra can be found 
in the books \cite{Frank-book,Fujishige-book}.  
For example, if $p$ is integer-valued, then the base-polyhedron $B'(p)$ 
is an integral polyhedron.  
Furthermore, the face of a base-polyhedron, its translation by a vector, 
its intersection with a box are also base-polyhedra, 
as well as the direct sum and the Minkowski sum of base-polyhedra.  
For an arbitrary chain ${\cal C}$ of subsets of $S$, the polyhedron
$\{x\in B'(p):  \ \widetilde x(Z) = p(Z)$ for each $Z\in {\cal C}\}$
is a face of $B'(p)$ and every face of $B'(p)$ arises in this way.

In applications it is fundamental that crossing supermodular (or
submodular) functions also define (possibly empty) base-polyhedra.
The non-emptiness in this case was characterized by an elegant theorem
of Fujishige \cite{Fujishige84e}.  
For a general overview of base-polyhedra, see the book \cite{Frank-book}.  
More recent results \cite{FrankJ67,FrankJ68} indicate that base-polyhedra 
defined by even weaker functions also show up in several applications.  
It is essential to see that formulating and proving results for
base-polyhedra is significantly easier when the bounding set-function
defining $B$ is fully supermodular (or submodular) but in applications
it is typical that $B$ is described by a weaker function. 
For such situations it is a typical task to modify appropriately the algorithm
developed for base-polyhedra given by fully supermodular functions.

Let $D=(V,A)$ be a loopless digraph endowed with integer-valued functions 
$f:A\rightarrow \underline{\bf Z}$ and 
$g:A\rightarrow \overline{\bf Z}$ 
for which $f\leq g$.  
Here $f$ and $g$ are serving as lower and upper bound functions, respectively.  
An edge $e$ is called {\bf tight} if $f(e)=g(e)$.  
The polyhedron $T(f,g):=\{x:  f\leq x\leq g\}$
is called a {\bf box}.

Let $\varrho_{D}(Z) = \varrho_{A}(Z)$ denote 
the number of edges of $D$ entering $Z\subseteq V$ 
while $\delta_{D}(Z) = \delta_{A}(Z)$ is 
the number of edges of $D$ leaving $Z$.  
For a function $x:A\rightarrow {\bf R}$,
$\varrho_{x}(v):= \sum [x(uv):uv\in A]$ and 
$\delta_{x}(v):= \sum [x(vu):  vu\in A]$.  
We call $x$ {\bf feasible} if $f\leq x\leq g$.  
It is a simple property that $f\leq g$ implies that 
$\varrho_{g}-\delta_{f}$ is a fully submodular function.  
Define function 
$\psi_{x}:V\rightarrow {\bf R}$ by 
$\psi_{x}(v)=\varrho_{x}(v)-\delta_{x}(v)$.  
The function $\psi_{x}$ on $V$ is sometime called the 
{\bf net in-flow} of $x$.  
Define the set-function $\Psi_{x}$ by 
$\Psi_{x}(Z) :=\varrho_{x}(Z)-\delta_{x}(Z)$ 
for $Z\subseteq V$.  
Clearly, $\Psi_{x}(\emptyset )= \Psi_{x}(V)=0$.

Let $b$ be a crossing submodular function.  
A function 
$x:A\rightarrow {\bf R}$ is called a {\bf submodular flow} 
if $\Psi_{x} \leq b$.  
When $f\leq x\leq g$, we speak of 
an {\bf $(f,g)$-bounded} or {\bf feasible} submodular flow.  
If there is a submodular flow $x$, then 
$\Psi_{x}(V)=0$ implies that $b(V)\geq 0$.  
If $b(V)>0$, then $b(V)$ can be reduced to $0$ 
since this affects neither the (crossing) submodularity of $b$ 
nor the requirement $\Psi_{x} \leq b$.  
Therefore we shall assume throughout that $b(V)=0$, 
that is, the base-polyhedron $B(b)$ is a $0$-base-polyhedron 
(with zero component-sum of each vector in $B$).  
The set 
\[
Q = \{x\in {\bf R}\sp A:  f\leq x\leq g, \Psi_{x}\leq b\}
\]
\noindent 
of $(f,g)$-feasible submodular flows is called a {\bf submodular flow} polyhedron.

It is immediate from the definitions that $x$ is a submodular flow
precisely if its net in-flow vector $\psi_{x}$ belongs to the
$0$-base-polyhedron $B$ described by $b$.  
Since a base-polyhedron can
also be defined with the help of a supermodular function, it follows
that submodular flows can be defined with supermodular functions, as
well, namely, if $p$ is a crossing supermodular function, then the
polyhedron 
\begin{equation} 
Q := \{x\in {\bf R}\sp A:  f\leq x\leq g, \Psi_{x}\geq p\} 
\label{(pbaseflow)} 
\end{equation} 
\noindent 
is also a submodular flow polyhedron.
Therefore the role of submodular and supermodular functions in
defining a submodular flow is completely symmetric, and hence a
submodular flow could also be called a supermodular flow.  
This is why we suggested in Section~\ref{main-results} the term base-flow
rather than submodular flow.  
Actually, the point is that, in the definition of a base-flow, 
it is the $0$-base-polyhedron $B$ that plays
the essential role and not the way how $B$ is given.  
Therefore a function (or vector) $x:A\rightarrow {\bf R}$ is a base-flow 
if $\psi_{x}\in B$
while a base-flow polyhedron 
\begin{equation} 
Q(f,g;B):=\{x\in {\bf R}\sp A:  \ f\leq x\leq g,\ \psi_{x} \in B\} 
\label{(QfgB)} 
\end{equation}
\noindent 
is the set of feasible base-flows.  
When the $0$-base-polyhedron $B=B'(p)$ is defined 
by a supermodular function $p$ (with $p(V)=0$), we speak of a $p$-base-flow.

A fundamental result of Edmonds and Giles \cite{Edmonds-Giles} states that, 
for a crossing submodular function $b$, the linear system
\[
 \{ f\leq x\leq g, \  \varrho_{x}(Z) -\delta_{x}(Z) \leq b(Z) 
 \ \ \hbox{for every}\ \ Z\subseteq V \} 
\]
\noindent 
describing a base-flow polyhedron is totally dual integral (TDI), 
implying that $Q$ is an integral polyhedron whenever $f$, $g$, $b$ are integer-valued.

An important feature of base-flows is that, given 
a subset $F\subseteq A$, 
the projection of a base-flow polyhedron $Q$ to ${\bf R} \sp F$
(that is, the restriction to $F$) is itself a base-flow polyhedron.
The face of a base-flow polyhedron is also a base-flow polyhedron.

When $B$ consists of a single element $m$ (with $\widetilde m(V)=0$),
a base-flow is a modular flow or mod-flow, 
which is a standard circulation when $m\equiv 0$.  
Decreasingly minimal integer-valued
mod-flows were investigated in \cite{Frank-Murota.C}.  
Note that the projection of a mod-flow polyhedron is 
a base-flow polyhedron but typically not a mod-flow polyhedron.

Frank \cite{FrankJ8,Frank-book} provided a necessary and sufficient
condition for the non-emptiness of a base-flow polyhedron.  
For general crossing supermodular (or submodular) functions this condition
is rather complicated but when $p$ is fully supermodular, the
formulation (and the proof of its necessity) is pretty straightforward.

\begin{theorem} 
\label{Q-non-empty} 
When $p$ is fully supermodular, the
$p$-base-flow polyhedron $Q$ defined in \eqref{(pbaseflow)} is non-empty 
if and only if \ 
$\varrho_{g}-\delta_{f}\geq p$, that is, 
\begin{equation}
\varrho_{g}(X)-\delta_{f}(X)\geq p(X) \quad
\hbox{\rm for every subset $X\subseteq V$}. 
\label{(baseflow-nonempty)} 
\end{equation} 
\end{theorem}

\par
Note that in the special case of circulations (when $p\equiv 0$), 
we are back at Hoffman's circulation theorem.

\section{$L$-upper-minimal base-flows} 
\label{flow3}

Let $D$, $f$, $g$, $B$, and $Q=Q(f,g;B)$ be the same as in Theorem~\ref{MAIN}, 
but in this section we do not use $F$.  
Let $L$ be a subset of $A$ for which $f$ and $g$ are finite-valued on $L$ 
(that is, $f(e)$ may be $-\infty $ and $g(e)$ may be $+\infty $ only if $e\in A-L$).  
We say that $z\in \odotZ{Q}$ is {\bf $L$-upper-minimal} or that
$z$ is an {\bf $L$-upper-minimizer} if the number of $g$-saturated edges in $L$ 
is as small as possible, where an edge $e\in L$ is called
{\bf $g$-saturated} if $z(e)=g(e)$.  
In this section, we are interested in characterizing 
the $L$-upper-minimizer elements of $\odotZ{Q}$.

For a chain ${\cal C}$ of subsets of $V$, we call an edge $e$ 
{\bf ${\cal C}$-entering}\ (respectively, 
{\bf ${\cal C}$-leaving}) 
if $e$ enters (resp., leaves) a member of ${\cal C}$, 
and $e$ is {\bf ${\cal C}$-neutral} 
when $e$ neither enters nor leaves any member of ${\cal C}$.
For a subset $I$ of edges, let $\varrho_{I}(C)$ \ ($C\subseteq V$)
denote the number of edges in $I$ entering $C$, and let 
$\varrho_{I}({\cal C})$ denote the number of ${\cal C}$-entering edges in $I$.
When this number is positive, we say that $I$ {\bf enters} ${\cal C}$.

One of the goals of this section is to prove the following
characterization of $L$-upper minimizer base-flows.  This will serve
as a main tool in proving Theorem~\ref{MAIN}.

\begin{theorem} 
\label{Lupmin} 
Let $Q=Q(f,g;B)$ be a non-empty integral
base-flow polyhedron and $L$ a subset of edges on which both $f$ and
$g$ are finite-valued.  There is a face $B_{L}$ of $B$ and there are bounds 
$f_{L}:A\rightarrow \underline{\bf Z}$ and $g_{L}:A\rightarrow \overline{\bf Z}$ 
with $f\leq f_{L}\leq g_{L}\leq g$ such that an element $z\in \odotZ{Q}$ 
is an $L$-upper minimizer if and only if $z\in \odotZ{Q_{L}}$ 
where $Q_{L}= Q(f_{L},g_{L};B_{L})$.  
\end{theorem}

The proof will be prepared in this section and completed in Section~\ref{minmaxL}.

\medskip 
If there is an edge $e\in L$ with $f(e)=g(e)$, then $x(e)=g(e)$ for each $x\in Q$.  
Therefore a base-flow $x$ is $g$-saturated on $L$ 
if and only if 
it is $g$-saturated on $L':=L-e$,
and hence $x$ is $L$-upper minimal precisely if it is $L'$-upper minimal.  
Therefore it suffices to prove the theorem for $L'$, and
hence we can assume that $L$ contains no tight edges, or in other words, 
\begin{equation} 
-\infty <f(e)<g(e)<+\infty \quad \hbox{for every edge}\ \ \ e\in L. 
\label{(Ldef)} 
\end{equation}

\par
We shall also show how the chain determining the face $B_{L}$ 
and the bounds $(f_{L},g_{L})$ occurring in the theorem can be computed in
polynomial time in the case when a subroutine is available to compute
$p(Z)$ for any given subset $Z\subseteq S$, where $p$ denotes 
the unique fully supermodular function $p$ for which $B=B'(p)$.

\subsection{Lower bound for the number of $g$-saturated edges}

Let $p$ denote the unique fully supermodular function defining $B$, 
that is, $B=B'(p)$.  
We say that a chain $\cal C$ of subsets of $V$ is 
{\bf feasible} 
if $\sum [\varrho_{g}(C) - \delta_{f}(C) -p(C): C \in {\cal C}] < +\infty$. 
In particular, this implies
that $g(e)$ is finite for each edge $e$ entering $\cal C$, 
$f(e)$ is finite for each edge leaving $\cal C$, 
and $p(C)$ is finite for each member $C$ of ${\cal C}$.

Our first goal is to show how a feasible chain provides 
a lower bound for the number of $g$-saturated edges in $L$.  
Let $B_{\cal C}$ be the face of $B$ defined by $\cal C$.  
Define a bounding pair $(f_{\cal C}(e),g_{\cal C}(e))$ for each edge $e\in A$, 
as follows.  
For $e\in L$, let
\begin{equation} 
(f_{\cal C}(e),g_{\cal C}(e)):= 
\begin{cases} 
(g(e),g(e)) & \hbox{if $e$ enters at least two members of ${\cal C}$},  \cr
(g(e)-1,g(e)) & \hbox{if $e$ enters exactly one member of ${\cal C}$},  \cr 
(f(e),f(e)) & \hbox{if $e$ is ${\cal C}$-leaving}, \cr
(f(e),g(e)-1) & \hbox{if $e$ is ${\cal C}$-neutral}. 
\end{cases}
\label{(f+g-def1)} 
\end{equation}
\noindent 
For $e\in A-L$, let 
\begin{equation} 
(f_{\cal C}(e),g_{\cal C}(e)):=
\begin{cases} 
(g(e),g(e)) & \hbox{if $e$ is ${\cal C}$-entering},  \cr
(f(e),f(e)) & \hbox{if $e$ is ${\cal C}$-leaving}, \cr 
(f(e),g(e)) & \hbox{if $e$ is ${\cal C}$-neutral}. 
\end{cases} 
\label{(f+g-def2)}
\end{equation}
\noindent 
Note that the feasibility of chain ${\cal C}$ implies that
$f_{\cal C}(e)$ is finite for each edge $e$ leaving $\cal C$ and
$g_{\cal C}(e)$ is finite for each edge $e$ entering $\cal C$.  
It follows from this definition that $f\leq f_{\cal C}\leq g_{\cal C}\leq g$.  
Note that these data define a base-flow polyhedron
$Q(f_{\cal C},g_{\cal C};B_{\cal C})$ included in $Q(f,g;B)$.

Consider the following optimality criteria:
\begin{equation} 
\begin{cases} 
\hbox{(O1)} \quad x(e)=f(e) 
& \ \  \hbox{if \quad $e\in A$ is ${\cal C}$-leaving}, \cr 
\hbox{(O2)} \quad x(e)=g(e) 
& \ \ \hbox{if \quad $e\in A-L$ is ${\cal C}$-entering}, \cr
\hbox{(O3)} \quad g(e)-1\leq x(e)\leq g(e)
& \ \ \hbox{if \quad $e\in L$ enters exactly one member of ${\cal C}$}, \cr 
\hbox{(O4)} \quad x(e)=g(e) 
& \ \ \hbox{if \quad $e\in L$ enters at least two members of ${\cal C}$}, \cr 
\hbox{(O5)} \quad f(e)\leq x(e)\leq g(e)-1 
& \ \ \hbox{if \quad $e\in L$ is ${\cal C}$-neutral}, \cr
\hbox{(O6)} \quad \varrho_{x}(C) - \delta_{x}(C) = p(C) 
& \ \ \hbox{if \quad $C\in {\cal C}$}. 
\end{cases} 
\label{(goptkrit)}
\end{equation}

An easy case-checking immediately shows the following.

\begin{claim}
\label{optkritekv} 
Let $x\in Q$.  
The union of the first five optimality criteria in \eqref{(goptkrit)} 
is equivalent to the requirement $f_{\cal C}\leq x\leq g_{\cal C}$.  
Criterion {\rm (O6)} is
equivalent to stating that $\psi_{x}\in B_{\cal C}$.  
\FBOX 
\end{claim}

\begin{lemma}  \label{est-crit} 
Let $Q=Q(f,g;B)$ be a non-empty base-flow polyhedron and 
let $p$ denote the unique supermodular function
defining the $0$-base-polyhedron $B$ (that is, $B=B'(p))$.  
Let $L$ be a subset of $A$ meeting \eqref{(Ldef)}.  
Let $x$ be an arbitrary element of $Q$ and let 
\begin{equation} 
X:=\{e\in L:  x(e)=g(e)\}.  
\label{(Xdef)}
\end{equation} 
\noindent
Let ${\cal C}$ be a feasible chain of subsets of $V$.  
Then 
\begin{equation} 
 | X | \geq \varrho_{L}({\cal C}) 
\ - \sum [\varrho_{g}(C)- \delta_{f}(C) - p(C):  \ C\in {\cal C}].  
\label{(ert-crit)} 
\end{equation}
\noindent 
Moreover, \eqref{(ert-crit)} is met by equality 
if and only if 
the optimality criteria hold in \eqref{(goptkrit)}, or, 
equivalently, $x\in Q(f_{\cal C},g_{\cal C};B_{\cal C})$.  
\end{lemma}

\Proof 
Since $\cal C$ is feasible, 
$\varrho_{g}(C)<+\infty$, $\delta _f(C)> -\infty$, 
and $p(C)>-\infty$ for each $C\in {\cal C}$.  
Observe that, for any set $I$ of edges, 
\begin{equation} 
\sum_{C\in {\cal C}}\varrho_{I}(C) \geq \varrho_{I}({\cal C}), 
\label{(Iineq)} 
\end{equation} 
(where $\varrho_{I}({\cal C})$ denotes 
the number of ${\cal C}$-entering elements of $I$)
and 
\begin{equation} 
\sum_{C\in {\cal C}} \varrho_{I}(C) 
= \varrho_{I}({\cal C}) \ \ \Leftrightarrow \ \ 
\hbox{each edge in $I$ enters at most one member of ${\cal C}$}. 
\label{(Ieq)} 
\end{equation}

By applying \eqref{(Iineq)} to $I:=L-X$, the assumption $x\in Q$ implies
\begin{align} 
\sum_{C\in {\cal C}} p(C) 
& \leq \sum_{C\in {\cal C}} [\varrho_{x}(C) - \delta_{x}(C)] 
\nonumber \\
& \leq \sum_{C\in {\cal C}} [\varrho_{g}(C) -\delta_{f}(C) -\varrho_{L-X}(C)] 
\nonumber \\
&\leq \sum_{C\in {\cal C}} [\varrho_{g}(C) - \delta_{f}(C)] - \varrho_{L-X}({\cal C}) 
\nonumber \\
&
= \sum_{C\in {\cal C}} [\varrho_{g}(C) - \delta_{f}(C)] -
   \varrho_{L}({\cal C}) + \varrho_X({\cal C}), 
\label{(roL)} 
\end{align}
\noindent from which
\begin{equation} 
| X | \geq \varrho_X({\cal C}) \geq \varrho_{L}({\cal C})
 + \sum_{C\in {\cal C}} [p(C) - \varrho_{g}(C) + \delta_{f}(C)],
\label{(vegso)} 
\end{equation} 
and hence \eqref{(ert-crit)} follows.  

\medskip\medskip

To see the second part of the lemma, suppose first that equality holds
in \eqref{(ert-crit)}.
We want to prove that the optimality criteria
hold in \eqref{(goptkrit)}.  
Equality in \eqref{(ert-crit)} implies that both inequalities in \eqref{(vegso)} are met with equalities.
Therefore 
\begin{equation} 
 | X | = \varrho_X({\cal C}) 
\label{(0.eq)}, 
\end{equation}
and each inequality in \eqref{(roL)} is met with equalities, that is,
\begin{align} 
\sum_{C\in {\cal C}} p(C) 
&= \sum_{C\in {\cal C}} [\varrho_{x}(C)-\delta_{x}(C)], 
\label{(1.eq)} 
\\
\sum_{C\in {\cal C}}
[\varrho_{x}(C) - \delta_{x}(C)] 
 &= \sum_{C\in {\cal C}} [\varrho_{g}(C) -\varrho_{L-X}(C) - \delta_{f}(C)], 
\label{(2.eq)} 
\\
\sum_{C\in {\cal C}} \varrho_{L-X}(C) 
&= \varrho_{L-X}({\cal C}) .
\label{(3.eq)} 
\end{align}
\noindent 
Here \eqref{(1.eq)} holds precisely if 
$p(C) = \varrho_{x}(C)- \delta_{x}(C)$ 
for each $C\in {\cal C},$ which is exactly (O6).

Equality \eqref{(2.eq)} implies that $\delta_{x}(C) = \delta_{f}(C)$ 
for each $C$ and hence (O1) follows.  
\eqref{(2.eq)} also implies for each ${\cal C}$-entering edge $e\in A-L$ 
that $x(e)=g(e)$, that is, (O2) holds.  
Furthermore, \eqref{(2.eq)} implies for each ${\cal C}$-entering
edge $e\in L$ that $g(e)-1\leq x(e) \leq g(e)$ (namely, $x(e)=g(e)$
when $e\in X$, and $x(e)=g(e)-1$ when $e\in L-X$), that is, (O3) holds.

Equalities \eqref{(3.eq)} and \eqref{(Ieq)} imply that $e\in L-X$ 
entering ${\cal C}$  enters exactly one member of ${\cal C}$.  
Therefore, if an edge $e\in L$ enters at least two members of ${\cal C}$, 
then $e\in X$, that is, $x(e)=g(e)$, showing that (O4) holds.

Finally, \eqref{(0.eq)} means that every $g$-saturated edge 
in $L$ enters ${\cal C}$.  
Therefore, if an edge $e\in L$ does not enter ${\cal C}$, 
then $x(e)\leq g(e)-1$, showing that (O5) holds.

\medskip

To see the reverse implication, suppose that the element $x$ of $Q$
and the chain ${\cal C}$ meet the six optimality criteria in the lemma.
We have to show that \eqref{(ert-crit)} holds with equality.  
What we are going to prove is that both inequalities in \eqref{(vegso)} 
are met with equality.

\begin{claim} 
\label{XroX} 
$ | X | = \varrho_X({\cal C})$.  
\end{claim}

\Proof 
Suppose, indirectly, that $ | X | > \varrho_X({\cal C})$,
which means that there is an edge $e\in X$ which is not ${\cal C}$-entering.  
Then $x(e)=g(e)>f(e)$ and (O1) imply that $e$ is not ${\cal C}$-leaving either.  
That is, $e$ is ${\cal C}$-neutral and 
hence (O5) implies that $x(e)\leq g(e)-1$, a contradiction.  
\FBOX
 
\medskip

Therefore the first inequality in \eqref{(vegso)} is met indeed by equality.  
The second inequality in \eqref{(vegso)} holds with equality 
precisely if each of the three inequalities in \eqref{(roL)} holds with equality.  
Consider these three inequalities separately.

Criterion (O6) states that $\varrho_{x}(C) - \delta_{x}(C) = p(C)$ 
for each $C\in {\cal C}$, implying \eqref{(1.eq)}, 
which shows that the first inequality in \eqref{(roL)} holds with equality.

\begin{claim} 
\label{2.ineq} 
The second inequality in \eqref{(roL)} holds with equality, that is, 
\begin{equation} 
\sum_{C\in {\cal C}} [\varrho_{x}(C) - \delta_{x}(C)] 
= \sum_{C\in {\cal C}} [\varrho_{g}(C) - \delta_{f}(C) -\varrho_{L-X}(C)].  
\label{(2.ineq)} 
\end{equation} 
\end{claim}

\Proof 
Observe first that (O1) implies 
\begin{equation} 
\sum_{C\in {\cal C}} \delta_{x}(C)= \sum_{C\in {\cal C}} \delta_{f}(C). 
 \label{(deltaxf)}
\end{equation}

Let $C$ be a member of ${\cal C}$.  
For an edge $e\in A-L$ entering $C$, (O2) implies that $x(e)=g(e)$.  
Furthermore, the definition of $X$ shows that 
$x(e)=g(e)$ for $e\in X$.  
By integrating these observations, we get the following:
\begin{align*}
& \varrho_{x}(C) 
\\
&
= \sum [x(e):  e\in A-(L-X),
\ e \ \hbox{enters}\ \ C] + 
\sum [x(e):  e\in L-X, \ e \ \hbox{enters}\ \ C]
\\
& = \sum [g(e):  e\in A-(L-X), \ e \ \hbox{enters}\ \ C] + 
    \sum [g(e)-1:  e\in L-X, \ e \ \hbox{enters}\ \ C]
\\ &=
 \varrho_{g}(C) - \varrho_{L-X}(C).
\end{align*}
\noindent
Therefore 
\[
\sum_{C\in {\cal C}} \varrho_{x}(C) = \sum_{C\in {\cal C}} [\varrho_{g}(C)- \varrho_{L-X}(C)].
\] 
\noindent
This and \eqref{(deltaxf)} imply \eqref{(2.ineq)}.  
\FBOX

\begin{claim}
 \label{roL-X} 
$\sum [\varrho_{L-X}(C):  C\in {\cal C}] = \varrho_{L-X}({\cal C)}$.  
\end{claim}

\Proof 
The claim is equivalent to stating that every edge in $L-X$ entering ${\cal C}$ 
enters exactly one member of ${\cal C}$.  
But this is true since  
if an edge $e\in L-X$, indirectly,  enters at least two members of ${\cal C}$, 
then (O4) implies that $x(e)=g(e)$, that is,
$e\in X$, a contradiction.  
\FBOX 

\medskip

Claim~\ref{roL-X} immediately implies that the third inequality in
\eqref{(roL)} also holds with equality.  
This and Claim~\ref{XroX} imply that 
both inequalities in \eqref{(vegso)} are met with equalities
and hence we have equality in \eqref{(ert-crit)}, as well.  
This completes the proof of Lemma \ref{est-crit}.
\BB

\section{Min-max formula for $L$-upper-minimizers} 
\label{minmaxL}

In this section we provide a min-max theorem for the minimum number of
$g$-saturated $L$-edges of a member of $\odotZ{Q}$.  
In Lemma~\ref{est-crit}, we did not need the unique supermodular function $p$
defining the $0$-base-polyhedron $B$ in question.  
However, in formulating and proving the min-max theorem below, we shall rely on $p$.  
It will also be shown how the optima can be computed with the
help of a standard algorithm to compute a cheapest integer-valued
feasible base-flow.

\begin{theorem} 
\label{minL} 
Let $Q=Q(f,g;B)$ be a non-empty base-flow polyhedron and 
let $p$ denote the unique supermodular function
defining $0$-base-polyhedron $B$ (that is, $B=B'(p))$.  
Let $L$ be a subset of $A$ meeting \eqref{(Ldef)}.  
The minimum number of
$g$-saturated $L$-edges of an $(f,g)$-bounded integer-valued
$p$-base-flow is equal to 
\begin{equation} 
\max \{ \varrho_{L}({\cal C}) \ 
  - \sum_{C\in {\cal C}} [\varrho_{g}(C)- \delta_{f}(C) - p(C)]\},
\label{(minmaxL)} 
\end{equation} 
where the maximum is taken over all feasible chains $\cal C$ of subsets of $V$.
In particular, if the minimum is zero,
the maximum is attained at the empty chain.  
\end{theorem}

\Proof 
The first part of Lemma~\ref{est-crit} implies that $\min \geq \max$.

For proving the reverse direction, consider the primal optimization problem,
which is to find an integral $(f,g)$-bounded base-flow 
which saturates (with respect to $g$) a minimum number of elements of $L$.
We show that this is a minimum cost feasible base-flow problem on a modified digraph.  
To this end, we introduce a parallel copy $e'$ of each $e\in L$.  
Let $L'$ denote the set of new edges.  
Let 
\[ 
A_{1}:=A\cup L', \quad D':=(V,L'), \quad
D_{1}:=(V,A_{1}).
\]
\noindent 
Define $g_{L}$ on $A$ by $g_{L}:=g-\chi_{L}$, that is, 
we reduce $g(e)$ by $1$ for each $e\in L$.  
Since $L$ contains no tight edges, $g_{L}\geq f$ follows.  
Let $f_{1}$ and $g_{1}$ be bounding functions and
$c_{1}$ a cost function defined on $A_{1}$ as follows:
\begin{align*}
& f_{1}(e) :  = f(e), \ \ g_{1}(e):=g_{L}(e), \ \ \ c_{1}(e):=0  \quad \hbox{if} \ \ e\in A,  
\\
& f_{1}(e') :  = 0, \ \ g_{1}(e'):= 1, \ \ \ c_{1}(e'):=1 \quad \hbox{if} \ \ e'\in L'.  
\end{align*}

\begin{claim} 
\label{ekvi} 
The problem of finding an $(f,g)$-bounded integer-valued $p$-base-flow on $A$ 
admitting a minimum number of $g$-saturated $L$-edges 
is equivalent to finding a minimum $c_{1}$-cost $(f_{1},g_{1})$-bounded 
integer-valued $p$-base-flow on $A_{1}$.  
\end{claim}

\Proof 
First, let $x\in \odotZ{Q}$ and let 
$X:=\{e\in L:  x(e)=g(e)\}$ denote the set of $g$-saturated elements of $L$.  
Let $X':=\{e':  e\in X\}$ denote the subset of $L'$ corresponding to $X$.  
Define a $p$-base-flow $x_{1}$ on $A_{1}$ as follows: 
\begin{align*}
& x_{1}(e):=
\begin{cases} 
x(e) & \hbox{if $e\in A-X$}, \cr 
g(e)-1 & \hbox{if $e\in X$}, 
\end{cases}
\\ & 
x_{1}(e'):= 
\begin{cases}
 1 & \hbox{if $e'\in X'$}, \cr 
0 & \hbox{if $e'\in L'-X'$}.  
\end{cases} 
\end{align*}
\noindent 
Then
$x_{1}$ is an $(f_{1},g_{1})$-bounded $p$-base-flow on $A_{1}$ 
whose $c_{1}$-cost is $ | X\vert$.

Conversely, let $x_{1}$ be a minimum cost 
integer-valued $(f_{1},g_{1})$-bounded $p$-base-flow on $A_{1}$.  
Observe that if
$x_{1}(e')=1$ for some $e'\in L'$, then $x_{1}(e)=g_{1}(e)=g(e)-1$
where $e$ is the edge in $L$ corresponding to $e'$.  
Indeed, if we had
$x_{1}(e)\leq g(e)-2$, then the $p$-base-flow obtained from $x_{1}$ 
by adding $1$ to $x_{1}(e)$ and subtracting $1$ from $x_{1}(e')$ would be of smaller cost.  
It follows that the $p$-base-flow $x$ on $A$ defined by
\begin{equation} 
x(e):= \begin{cases} 
x_{1}(e)+x_{1}(e') & \hbox{if $e\in L$}, \cr 
x_{1}(e) & \hbox{if $e\in A-L$}
\end{cases}
\label{(xdef)} 
\end{equation}
is an $(f,g)$-bounded $p$-base-flow in $D$, 
for which the number of $g$-saturated $L$-edges is exactly the $c_{1}$-cost of $x_{1}$.  
\FBOX

\medskip

By Claim~\ref{ekvi}, we investigate the minimum of the $c_{1}$-cost of
integer-valued $(f_{1},g_{1})$-bounded $p$-base-flows on $A_{1}$.  
In order to describe the dual optimization problem, 
let $N$ be a $\{0,\pm 1\}$-matrix whose columns correspond to the elements of $A$ 
(the edge-set of $D$) while its rows correspond to the members of 
\[
{\cal P} :=\{Z\subseteq V:  \ p(Z)>-\infty \}.
\] 
\noindent 
The entry of $N$ corresponding to $Z\in {\cal P}$ and $e\in A$ is 
$+1$ if $e$ enters $Z$, 
$-1$ if $e$ leaves $Z$, and 0 otherwise.  
Note that, for a vector
$x\in A\rightarrow {\bf R}$,  
the requirement $\varrho_{x}(Z)-\delta_{x}(Z)\geq p(Z)$ 
for every $Z\subseteq V$ is equivalent to $Nx \geq p$.
Let $N'$ denote the matrix associated analogously with $D'$,  and let $N_{1}=[N,N']$.

The primal linear program is as follows:  
\begin{equation} 
\min \{c_{1}x_{1} :  \ N_{1}x_{1}\geq p, \ x_{1}\geq f_{1}, \ -x_{1} \geq -g_{1}\}.
\label{(primalflow)} 
\end{equation}
The dual linear program is as follows:  
\begin{equation} 
\max \{yp + z_{1}f_{1} - w_{1}g_{1}: 
 \ yN_{1} + z_{1} - w_{1} = c_{1}, \ y\geq 0, \ z_{1}\geq 0, \ w_{1}\geq 0\} 
\label{(dualflow)} 
\end{equation}
\noindent 
where $yp = \sum [y(Z)p(Z):  \ Z \in {\cal P}]$.  
Note that the components of $z_{1}=(z,z')$ 
correspond to the edges in $A$ and in $L'$, respectively, 
and the analogous statement holds for $w_{1}=(w,w')$.

A fundamental theorem of Edmonds and Giles \cite{Edmonds-Giles} states
that the linear system in \eqref{(primalflow)} is totally dual integral (TDI).  
In the present case, when each of $f,g,p,$ and $c_{1}$ is
integer-valued, the TDI-ness implies that both the primal and the dual
program have an integer-valued optimal solution.  
In addition, since $p$ is fully supermodular, 
the optimal $y\sp{*}$ in \eqref{(dualflow)}
can be chosen in such a way 
that the sets $C$ for which $y\sp{*}(C)$ is positive form a chain.  
Therefore we consider only those vectors $y$
occurring in the constraint of \eqref{(dualflow)} which are {\bf chained} 
in the sense that ${\cal C}=\{C\subseteq V:  y(C)>0\}$ is a chain.  
We refer to ${\cal C}$ as the {\bf support chain} of $y$.  
The dual solution $(y,z_{1},w_{1})$ is said to be {\bf simple} if $y$ is
chained and $\min \{z_{1}(e), w_{1}(e)\}=0$ for every edge $e\in A_{1}$.

\begin{claim} 
\label{vansimple} 
There exists an integer-valued optimal dual solution to \eqref{(dualflow)}
 which is simple.  
\end{claim}

\Proof 
We indicated already that there is a dual integral solution
$(y,z_{1},w_{1})$ where $y$ is chained.  
If both $z_{1}(e)$ and $w_{1}(e)$ are positive on an edge $e\in A_{1}$, 
then by reducing both
$z_{1}(e)$ and $w_{1}(e)$ by $\min \{z_{1}(e), w_{1}(e)\}$ 
we obtain another dual solution whose dual cost is larger by
$g_{1}(e)-f_{1}(e)\geq 0$ 
than the dual cost $yp + z_{1}f_{1} - w_{1}g_{1}$ of $(y,z_{1},w_{1})$.  
\FBOX

\medskip

Let $(y,z_{1},w_{1})$ be a simple (integer-valued) dual solution.
Since $z_{1}$ and $w_{1}$ are non-negative vectors, 
$y$ uniquely determines them.  
We describe explicitly how $z_{1}$ and $w_{1}$ can be expressed by $y$.  
Recall that chain ${\cal C}$ consist of subsets $Z\subseteq V$ 
for which $y(Z)$ is positive.  
Define a function $\pi_{y}:V\rightarrow {\bf Z}$ by
\begin{equation} 
\pi_{y}(v) := \sum [y(C) : v\in C\in {\cal C}].  
\label{(pidef)} 
\end{equation}
\noindent 
Also define $\Delta_{y}:A_{1}\rightarrow {\bf Z}$ as follows:  
For $a=uv\in A_{1}$, let 
\[ 
\Delta_{y}(a):= \pi_{y}(v)-\pi_{y}(u).
\] 
\\
\noindent 
Then $\Delta_{y}(a) >0$ (respectively, $\Delta_{y}(a) < 0$) 
if and only if
$a$ is ${\cal C}$-entering (resp., ${\cal C}$-leaving) 
and hence $\Delta_{y}(a)=0$ precisely when $a$ is ${\cal C}$-neutral.

\begin{claim} 
\label{pidiff} 
Let $a=uv$ be an edge of $D_{1}$.  Then 
\begin{equation} 
\Delta_{y}(a) + z_{1}(a) -w_{1}(a) = c_{1}(a).  
\label{(pidiff)} 
\end{equation} 
\end{claim}

\Proof 
It follows from the definition of $\Delta_{y}$ that 
\begin{align}
&\hbox{$\Delta_{y}(a) = \sum [y(C):  \ C\in {\cal C}$ is entered by $a]$} 
\quad \,
\hbox{when $a$ is ${\cal C}$-entering}, 
\label{(TenEnt)} 
\\
&\hbox{$\Delta_{y}(a) = -\sum [y(C):  \ C\in {\cal C}$ is left by $a]$} 
\qquad   
\hbox{when $a$ is ${\cal C}$-leaving}, 
\label{(TenLea)} 
\\
&\hbox{$\Delta_{y}(a)=0$} 
\qquad 
\hbox{when $a$ is ${\cal C}$-neutral}.
\label{(TenNeu)} 
\end{align}
\noindent 
These and the constraint $yN_{1} + z_{1} - w_{1} = c_{1}$ 
in \eqref{(dualflow)} imply \eqref{(pidiff)}.  
\FBOX

\begin{claim} 
\label{zw.A} 
Let $e=uv\in A$.  
\begin{align}
&\hbox{If $e$ is ${\cal C}$-neutral, then} \quad 
\begin{cases} 
z_{1}(e)=0, &  \cr 
w_{1}(e)=0, &
\end{cases} 
\\
&\hbox{if $e$ is ${\cal C}$-entering, then} \quad
\begin{cases} 
z_{1}(e)=0, & \ \cr 
w_{1}(e)= \Delta_{y}(e), 
\end{cases}
\\
&\hbox{if $e$ is ${\cal C}$-leaving, then} \quad
\begin{cases} 
z_{1}(e)= -\Delta_{y}(e) , \cr
w_{1}(e)=0.   
\end{cases}
\end{align}
\end{claim}

\Proof 
Since $c_{1}(e)=0$ and $(y,z_{1},w_{1})$ is simple, the claim follows
by applying \eqref{(pidiff)} to $e$ in place of $a$.  
\FBOX

\begin{claim} \label{zw.L'} Let $e'=uv\in L'$.  
\begin{align}
&\hbox{If $e'$ is ${\cal C}$-neutral, then} \quad 
\begin{cases} 
z_{1}(e')=1,   \cr 
w_{1}(e')=0,  
\end{cases} 
\\
&\hbox{if $e'$ is ${\cal C}$-entering, then} \quad
\begin{cases} 
z_{1}(e')=0, \cr 
w_{1}(e')= \Delta_{y}(e') -1, 
\end{cases} 
\\ 
&\hbox{if $e'$ is ${\cal C}$-leaving, then} \quad
\begin{cases} 
z_{1}(e')= -\Delta_{y}(e')+1, \cr
w_{1}(e')=0. 
\end{cases} 
\end{align}
\end{claim}

\Proof 
Since $c_{1}(e')=1$ and $(y,z_{1},w_{1})$ is simple, the claim follows 
by applying \eqref{(pidiff)} to $e'$ in place of $a$.  
\FBOX

\begin{claim} 
\label{z1f1} 
$z_{1}f_{1} = \sum [y(C)\delta_{f}(C):  C\in {\cal C}]$.  
\end{claim}

\Proof 
Since $f_{1}(e')=0$ for $e'\in L'$, we have 
\begin{align*}
z_{1}f_{1} 
&= \sum_{e\in A}z_{1}(e)f_{1}(e) + \sum_{e'\in L'} z_{1}(e')f_{1}(e')
 = \sum_{e\in A}z_{1}(e)f_{1}(e)
\\ &=
\sum [z_{1}(e)f(e):  \ e\in A, \ z_{1}(e)>0] = \sum
[ - \Delta_{y}(e) f(e):  \ e\in A \ \hbox{is ${\cal C}$-leaving$]$}
\\ &=
\sum [y(C)\delta_{f}(C):  C\in {\cal C}].
\end{align*}
\FBOX

\begin{claim} 
\label{w1g1} 
$w_{1}g_{1} = \sum [ y(C)\varrho_{g}(C):  C\in {\cal C}] - \varrho_{L}({\cal C})$.
\end{claim}

\Proof 
Observe first that 
\begin{align*}
& \sum [ \Delta_{y}(e) g_{1}(e):  \ e\in A \ \hbox{is ${\cal C}$-entering$]$}
\\ & 
= \sum [\Delta_{y}(e) g(e): e\in A-L \ \hbox{is ${\cal C}$-entering$]$}
  + \sum [\Delta_{y}(e) (g(e)-1):  e\in L \ \hbox{is ${\cal C}$-entering$]$} 
\\ &
= \sum [
\Delta_{y}(e) g(e):  \ e\in A \ \hbox{is ${\cal C}$-entering$]$}
 - \sum [ \Delta_{y}(e) :  \ e\in L \ \hbox{is ${\cal C}$-entering$]$}
\\ & 
= \sum_{C\in {\cal C}} y(C)\varrho_{g}(C)
  - \sum [ \Delta_{y}(e) : \ e\in L \ \hbox{is ${\cal C}$-entering$]$}.
\end{align*}
Second, 
\begin{align*}
& \sum [\Delta_{y}(e')-1:  \ e'\in L' \ \hbox{is ${\cal C}$-entering$]$}
\\ &=
 \sum [\Delta_{y}(e)-1:  \ e\in L \ \hbox{is ${\cal C}$-entering$]$}
\\ &=
\sum [\Delta_{y}(e):  \ e\in L \ \hbox{is
${\cal C}$-entering$]$}\ - \sum [1 :  e\in L \ \hbox{is ${\cal C}$-entering$]$} 
\\ &=
\sum [\Delta_{y}(e):  \ e\in L \ \hbox{is ${\cal C}$-entering$]$} 
    - \varrho_{L}({\cal C}).  
\end{align*}
Since $g_{1}(e')=1$ for $e'\in L'$, we have 
\begin{align*}
w_{1}g_{1} 
 &=
\sum_{e\in A}w_{1}(e)g_{1}(e) 
 + \sum_{e'\in L'}w_{1}(e')g_{1}(e') 
\\ &=
 \sum [w_{1}(e)g_{1}(e): e\in A, w_{1}(e)>0] + \sum_{e'\in L'}w_{1}(e')
\\ &=
\sum [ \Delta_{y}(e) g_{1}(e):  \ e\in A \ \hbox{is ${\cal C}$-entering$]$}
  + \sum [\Delta_{y}(e') -1:  \ e'\in L'\ \hbox{is ${\cal C}$-entering$]$} 
\\ &=
\sum_{C\in {\cal C}} y(C)\varrho_{g}(C)
  - \sum [ \Delta_{y}(e) :  \ e\in L \ \hbox{is ${\cal C}$-entering$]$}
\\ & \qquad \qquad 
 {}+ \sum [\Delta_{y}(e):  \ e\in L \ \hbox{is ${\cal C}$-entering$]$}
   - \varrho_{L}({\cal C}) 
\\ &=
\sum_{C\in {\cal C}} y(C)\varrho_{g}(C)- \varrho_{L}({\cal C}).
\end{align*}
\FBOX

\begin{claim} 
\label{dualcost} 
For a function $y:{\cal P}\rightarrow {\bf Z}_{+}$, 
let $z_{1}:A_{1}\rightarrow {\bf Z}_{+}$ and $w_{1}:A_{1}\rightarrow {\bf Z}_{+}$ 
be the vectors associated with $y$ in 
Claims {\rm \ref{zw.A}} and {\rm \ref{zw.L'}}.  
Then the cost of the dual solution $(y,z_{1},w_{1})$ in \eqref{(dualflow)} is as follows:  
\begin{equation} 
yp+ z_{1}f_{1} - w_{1}g_{1} 
= \varrho_{L}({\cal C}) - \sum_{C\in {\cal C}} y(C) [\varrho_{g}(C) -\delta_{f}(C)-p(C)].  
\label{(dualcost)} 
\end{equation} 
\end{claim}

\Proof 
By Claims \ref{z1f1} and \ref{w1g1}, we have 
\begin{align*}
& yp+ z_{1}f_{1} - w_{1}g_{1} 
\\ &=
\sum_{C\in {\cal C}} y(C)p(C) 
  + \sum_{C\in {\cal C}} y(C)\delta_{f}(C) 
  - \sum_{C\in {\cal C}} y(C)\varrho_{g}(C) + \varrho_{L}({\cal C)} 
\\ &=
\varrho_{L}({\cal C}) 
  - \sum_{C\in {\cal C}} y(C) [\varrho_{g}(C) -\delta_{f}(C) - p(C)].
\end{align*}
\FBOX

\begin{claim} 
\label{yreduce} 
Let $(y,z_{1},w_{1})$ be a simple integer-valued (possibly not optimal) solution 
to the dual problem \eqref{(dualflow)}.  
Let ${\cal C}$ denote the support chain of $y$,
and let $z_{1}$ and $w_{1}$ be the vectors associated with $y$ in 
Claims {\rm \ref{zw.A}} and {\rm \ref{zw.L'}}.  
Assume that $y(C_{1})\geq 2$ for a
member $C_{1}$ of ${\cal C}$.  Let $y'$ be the function arising from $y$
by reducing the value $y(C_{1})$ by $1$, and let $z'_{1}$ and $w'_{1}$ be
the vectors associated with $y'$ as described in 
Claims {\rm \ref{zw.A}} and {\rm \ref{zw.L'}}.  
Then 
\[ 
y'p+ z'_{1}f_{1} - w'_{1}g_{1} \geq yp+ z_{1}f_{1} - w_{1}g_{1}.
\] 
\end{claim}

\Proof 
By applying Claim~\ref{dualcost} to $(y,z_{1},w_{1})$ and to
$(y',z'_{1},w'_{1})$, we get  
\[ 
y'p+ z'_{1}f_{1} - w'_{1}g_{1} 
= yp+ z_{1}f_{1} - w_{1}g_{1} + [\varrho_{g}(C_{1}) -\delta_{f}(C_{1}) - p(C_{1})] 
\geq yp+ z_{1}f_{1} - w_{1}g_{1},
\] 
\noindent 
where the last inequality follows from \eqref{(baseflow-nonempty)}.  
\FBOX


\begin{proposition} 
\label{yis1} 
Let $(y,z_{1},w_{1})$ be a simple integer-valued optimal solution 
to the dual problem \eqref{(dualflow)} with support chain ${\cal C}_{L}$ of $y$, 
where $z_{1}$ and $w_{1}$ are the vectors associated with $y$ 
in Claims {\rm \ref{zw.A}} and {\rm \ref{zw.L'}}.  
Define the function $y_{L}:  {\cal P}\rightarrow \{0,1\}$ to be $1$ 
on the members of ${\cal C}_{L}$ and $0$ otherwise, 
and let $z'_{1}$ and $w'_{1}$ be the vectors 
associated with $y_{L}$ as described in  Claims {\rm \ref{zw.A}} and {\rm \ref{zw.L'}}.  
Then $(y_{L},z'_{1},w'_{1})$ is also a simple integer-valued optimal solution 
to the dual problem \eqref{(dualflow)}.  
\end{proposition} 

\Proof 
By repeated applications of Claim~\ref{yreduce}, we obtain that
$y_{L}p+ z'_{1}f_{1} - w'_{1}g_{1} \geq yp+ z_{1}f_{1} - w_{1}g_{1}$, 
but here we must have equality due to the hypothesis that 
$(y,z_{1},w_{1})$ is a dual optimum, 
and this equality shows that $(y_{L},z'_{1},w'_{1})$ is also a dual optimum.
\FBOX 

\medskip

Let $x\sp{*}_{1}$ be an integer-valued primal optimum in \eqref{(primalflow)}, 
that is, $x\sp{*}_{1}$ is a minimum $c_{1}$-cost
$(f_{1},g_{1})$-bounded $p$-base-flow on the edge-set $A_{1}$ of $D_{1}$.
Let $x\sp{*}$ be the $(f,g)$-bounded $p$-base-flow on $D$ defined in \eqref{(xdef)}.  
As noted in the proof of Claim~\ref{ekvi}, $x\sp{*}$
is an $L$-upper-minimizer, and the number of $g$-saturated $L$-edges 
of $x$ is $c_{1}x_{1}$.

Let $(y_{L},z_{1}',w_{1}')$ be a dual optimum in \eqref{(dualflow)} ensured
by Proposition~\ref{yis1}, where $y_{L}$ is $(0,1)$-valued with support
chain ${\cal C}_{L}$.  By Claim~\ref{dualcost}, 
the theorem of Edmonds and Giles
implies that the minimum number of $g$-saturated
$L$-edges of an $(f,g)$-bounded integer-valued $p$-base-flow is equal to 
\[ 
c_{1}x_{1}\sp{*} = y_{L}p+ z_{1}'f_{1} - w_{1}'g_{1} 
= \varrho_{L}({\cal C}_{L}) - \sum_{C\in {\cal C}} [\varrho_{g}(C) -\delta_{f}(C) -p(C)],
\]
completing the proof of Theorem~\ref{minL}.  
\BB

\medskip \medskip

\noindent 
{\bf Proof of Theorem~\ref{Lupmin}}.  \ 
Let ${\cal C}_{L}$ be an optimal chain in Theorem~\ref{minL}, 
and let $(f_{L},g_{L})$ be the pair of bounding functions associated with ${\cal C}_{L}$ 
in \eqref{(f+g-def1)} and \eqref{(f+g-def2)}
 (with ${\cal C}_{L}$ in place of ${\cal C}$), 
that is, 
$f_{L}:=f_{{\cal C}_{L}}$ and $g_{L}:=g_{{\cal C}_{L}}$.  
Let $B_{L}$ \ ($:=B_{{\cal C}_{L}})$ be the face of $B$ 
determined by ${\cal C}_{L}$, and let $Q_{L}:=Q(f_{L},g_{L};B_{L})$.
Theorem~\ref{minL} and Lemma~\ref{est-crit} imply that 
an element $x$ of $\odotZ{Q}$ is $L$-upper minimal 
if and only if 
$x\in \odotZ{Q_{L}}$, as required for Theorem~\ref{Lupmin}.  
\FBOX

\medskip

\begin{remark} \rm
The proofs in this section and Section~\ref{flow3} are a
bit technical and lengthy but the underlying ideas are pretty
standard.  Our objective of this detailed description is to make
easier to check the proofs.  
We note that the second part of Lemma~\ref{est-crit} can also 
be proved on the basis of the complementary slackness condition 
for the dual linear programs \eqref{(primalflow)}
and \eqref{(dualflow)} in the proof of Theorem~\ref{minL}.  
We also remark that the approach of introducing parallel edges in the proof of
Theorem~\ref{minL} is a well-known technique and is equivalent to
considering a convex cost which is equal to $1$ if $x(e)=g(e)$ and $0$
if $f(e)\leq x(e)\leq g(e)-1$.  
$\bullet $ 
\end{remark}


\noindent 
{\bf Algorithmic aspects}.  
First we outline how the proof above for Theorem~\ref{Lupmin}
gives rise to an algorithm for the case 
when the integral base-polyhedron $B$ is given by 
its unique fully supermodular defining function $p$, that is, $B=B'(p)$.  
After this description, we shall indicate when and how
this algorithm can be implemented in the case when $B$ is described 
in some other ways, for example, by a crossing supermodular function.

Since the cost-function $c_{1}$ is $\{ 0,1 \}$-valued, 
the algorithm of Cunningham and Frank \cite{FrankJ15} computes
in strongly polynomial time these primal and dual optima, provided a
subroutine is available to minimize a submodular function.  A simpler
and more efficient submodular flow algorithm, based on a push-relabel
approach, appeared in \cite{FrankJ62}.  
These algorithms actually work for the more general setting when the
bounding set-function is only crossing (sub- or) supermodular.  
They output an optimal dual solution $y$ whose support family 
${\cal F}_y:= \{ Z\subset V :  y(Z)>0\}$ is cross-free.  
But in the present case when $p$ is fully supermodular, 
it is straightforward to modify the algorithms 
so as to output an optimal dual solution $y$ such that 
the sets on which $y$ is positive form a chain ${\cal C}_{L}$ of 
non-empty proper subsets of $V$.  
Recall that, for determining the base-flow polyhedron 
$Q_{L}=Q(f_{L},g_{L};B_{L})$, all what we need is 
this chained optimal dual solution $y$ since this determines its support chain
${\cal C}_{L}$, which in turn gives rise to $(f_{L},g_{L})$, 
as given in \eqref{(f+g-def1)} and \eqref{(f+g-def2)}.

Consider now the case when the non-empty base-polyhedron $B$ is given
in some other way (for example, by a crossing supermodular function).
The above algorithmic approach concerning fully supermodular functions
can be applied in this case as well, since, 
by a basic property of base-polyhedra, 
$B$ determines its unique fully supermodular bounding function $p$, as follows:
\begin{equation} 
p(Z)=\min \{\widetilde x(Z):  x\in B\} .
\label{(moho)} 
\end{equation} 
\noindent 
Therefore the only requirement is that a subroutine be available 
to compute $p(Z)$ for any input set $Z\subseteq V$.  
Note that a minimization problem over a base-polyhedron given 
in \eqref{(moho)} can be solved with the help of an extension 
of Edmonds' greedy algorithm concerning polymatroids.  
See, for example, Section~14.5 in \cite{Frank-book} 
along with the geometric view of the greedy algorithm on Page 488, 
which explains how the greedy algorithm can be implemented 
for base-polyhedra given implicitly.  
All what we need is an oracle for minimizing a submodular function.

We also remark that the submodular flow (= base-flow) algorithms 
in the literature work with crossing submodular functions.  
But these algorithms can be extended to the case when the $0$-base-polyhedron $B$ 
defining the base-flow polyhedron is given in an implicit form.
What is really needed is a subroutine that is capable to find 
an integral element of $B$, and to decide for an element $z$ of $\odotZ{B}$ 
and a pair $\{s,t\}$ of elements of the ground-set whether 
$z':= z -\chi_{s} + \chi_{t}$ is in $B$ or not.

\section{Description of $F$-dec-min base-flows} \label{decmin}

Let $D=(V,A)$, $F$, $(f,g)$, and $B$ be the same as in Theorem~\ref{MAIN}.  
In this section, we prove this theorem by describing the set 
of $F$-dec-min integral elements of a non-empty integral base-flow polyhedron $Q=Q(f,g;B)$.  
When $F=\emptyset $, each element of $\odotZ{Q}$ is dec-min.  
In this case, the theorem is trivially true, so we assume henceforth that $F$ is non-empty.

A natural reduction step consists of removing a tight edge $e$ from $F$.  
This simply means that we replace $F$ by $F':=F-e$ 
(but keep $e$ in the digraph itself).  
Obviously, a member $z$ of $\odotZ{Q}$ is $F$-dec-min 
if and only if 
$z$ is $F'$-dec-min.  
Therefore, we may always assume that $F$ contains no tight edges, that is, 
\begin{equation} 
f(e) < g(e) \quad \hbox{for every $e\in F$}.
\label{(notight)} 
\end{equation} 
Removing tight edges from $F$ will be used not only 
at the starting step of the proof (and the algorithm) 
but it is a basic tool in later phases, as well, 
when the bounds $(f,g)$ are tightened and new tight edges arise in $F$.  
It is this reduction step that makes the current $F$ smaller
and smaller (see, Theorem~\ref{reduction}).

Moreover, we assume that both $f$ and $g$ are finite-valued on $F$, that is, 
\begin{equation} 
-\infty <f(e)<g(e)<+\infty \quad \hbox{for every $e\in F$}, 
\label{(fgveges)} 
\end{equation} 
\noindent 
which ensures the existence of an $F$-dec-min element of $Q$.  
When \eqref{(fgveges)} is not assumed, it
is possible that $\odotZ{Q}$ has no $F$-dec-min element at all.  
We shall characterize this situation in Section~\ref{vandecmin}.

\subsection{Pre-decreasing minimality on $F$} 
\label{predecminF}

Let $p$ denote the unique fully supermodular function defining 
the $0$-base-polyhedron $B$, that is, $B=B'(p)$.  
Let $\beta =\beta_F$ denote the smallest integer 
for which $\odotZ{Q}$ has an element $z$
that is {\bf $\beta $-covered on $F$} 
(meaning that $z(e)\leq \beta$ for every edge $e\in F$), 
that is, 
\begin{equation} 
\beta_F = \min \{ \max \{z(e): e\in F \}:  z\in \odotZ{Q}(f,g;B)\}.  
\label{(betaF)} 
\end{equation} 
\noindent 
In the next section, we shall work out an algorithm to compute $\beta_F$ 
in strongly polynomial time.  
Since we are interested in $F$-dec-min members of $\odotZ{Q}$, 
we may assume that the largest $g$-value of the edges in $F$ 
is this $\beta $. 
Let 
\begin{equation} 
L_\beta :=\{e\in F:  g(e)=\beta \}.  
\label{(Lbetadef)} 
\end{equation} 
\noindent 
By the definition of $\beta $, $L_\beta $ is non-empty.  
Now Condition \eqref{(baseflow-nonempty)} holds but,
since $F$ contains no tight edges and since $\beta $ is minimal, 
after decreasing the $g$-value of the elements of $L_\beta $ 
from $\beta $ to $\beta -1$, the resulting function 
$g\sp{-}:=g-\chi_{L_\beta }$
violates \eqref{(baseflow-nonempty)}, that is, 
$Q(f,g\sp{-};B)=\emptyset$. 
Summing up, we shall rely on the following notation and assumptions.

\begin{equation} 
\begin{cases}
 & \hbox{$F$ is non-empty and contains no
$(f,g)$-tight edges,}\ \cr 
& \hbox{$\beta :=\max \{g(e):  e\in F\}$, }\ \cr 
& \hbox{$L_\beta :=\{e\in F:  g(e) =\beta \}$, }\ \cr 
&\hbox{$g\sp{-}:=g-\chi_{L_\beta }$,}\ \cr 
& \hbox{$\odotZ{Q} = \odotZ{Q}(f,g;B)$ is non-empty, }\ \cr 
& \hbox{$\odotZ {Q}(f,g\sp{-};B)$ is empty.  }\ 
\end{cases} 
\label{(hypo)} 
\end{equation}

\medskip

As a preparation for deriving the main result Theorem~\ref{MAIN}, 
we need the following relaxation of decreasing minimality.  
We call a member $z$ of $\odotZ{Q}$ {\bf pre-decreasingly minimal} 
({\bf pre-dec-min}, for short) {\bf on} $F$ 
if the number $\mu$ of edges $e$ in $L_\beta $ 
with $z(e) = \beta $ is as small as possible.
Obviously, $z$ is pre-dec-min on $F$ precisely if $z$ is $L_\beta$-upper minimizer.  
It is also straightforward that if $z$ is $F$-dec-min, then $z$ is pre-dec-min on $F$.

Apply Theorem~\ref{Lupmin} to this $L:=L_\beta $, and consider 
the base-flow polyhedron 
\[
Q_{L}=(f_{L},g_{L};B_{L})
\] 
ensured by the theorem.
Recall that the face $B_{L}$ of $B$ was defined by a 
feasible chain ${\cal C}_{L}$
while $(f_{L},g_{L})$ was defined in \eqref{(f+g-def1)} and \eqref{(f+g-def2)}.  
In the present special case of $L=L_\beta $, 
the definition of $(f_{L},g_{L})$ in \eqref{(f+g-def1)} and \eqref{(f+g-def2)}
specializes as follows.  For $e\in L$, let
\begin{equation} 
(f_{L}(e),g_{L}(e)):= \begin{cases}
 (\beta ,\beta ) 
& \hbox{if $e$ enters at least two members of ${\cal C}_{L}$}, \cr 
(\beta -1,\beta ) &
\hbox{if $e$ enters exactly one member of ${\cal C}_{L}$}, \cr
(f(e),f(e)) 
& \hbox{if $e$ is ${\cal C}_{L}$-leaving}, \cr 
(f(e),\beta -1) 
& \hbox{if $e$ is ${\cal C}_{L}$-neutral}. 
\end{cases}
\label{(f'g'def1)} 
\end{equation}
\noindent 
For $e\in A-L$, let 
\begin{equation} 
(f_{L}(e),g_{L}(e)):= \begin{cases}
(g(e),g(e)) & \hbox{if $e$ is ${\cal C}_{L}$-entering}, \cr 
(f(e),f(e)) & \hbox{if $e$ is ${\cal C}_{L}$-leaving},  \cr 
(f(e),g(e)) & \hbox{if $e$ is ${\cal C}_{L}$-neutral}. 
\end{cases}
 \label{(f'g'def2)} 
\end{equation}

The optimality criteria \eqref{(goptkrit)}, when applied to 
$z\in \odotZ{Q}$ in place of $x$, are as follows:
\begin{equation} 
\begin{cases} 
\hbox{(O1)}\ \ \ \ \ z(e)=f(e) 
 & \ \ \hbox{if \quad $e\in A$ is ${\cal C}_{L}$-leaving}, \cr 
\hbox{(O2)}\ \ \ \ \ z(e)=g(e)
& \ \ \hbox{if \quad $e\in A-L$ is ${\cal C}_{L}$-entering}, \cr
\hbox{(O3)}\ \ \ \ \ \beta -1\leq z(e)\leq \beta 
& \ \ \hbox{if \quad $e\in L$ enters exactly one member of ${\cal C}_{L}$}, \cr 
\hbox{(O4)}\ \ \ \ \ z(e)=\beta 
& \ \ \hbox{if \quad $e\in L$ enters at least two members of ${\cal C}_{L}$}, \cr 
\hbox{(O5)}\ \ \ \ \ f(e)\leq z(e)\leq \beta -1 
& \ \ \hbox{if \quad $e\in L$ is ${\cal C}_{L}$-neutral}, \cr
\hbox{(O6)}\ \ \ \ \ \varrho_{z}(Z) - \delta_{z}(Z) = p(Z) 
& \ \ \hbox{if \quad $Z\in {\cal C}_{L}$}. 
\end{cases}
\label{(betaoptkrit)}
\end{equation}

By Lemma~\ref{est-crit} and Theorem~\ref{Lupmin}, we have the following.

\begin{claim} 
\label{Fdecmin} 
For an element $z\in \odotZ{Q}$, the following properties are equivalent:

\noindent 
{\rm (A)} \ $z$ is pre-dec-min on $F$,

\noindent 
{\rm (B)} \  \eqref{(betaoptkrit)} holds,

\noindent 
{\rm (C)} \ $z\in \odotZ{Q_{L}}$.  
\FBOX 
\end{claim}

\begin{claim} 
\label{f'g'} \ 
An element $z$ of $\odotZ{Q}$ is $F$-dec-min if
and only if $z$ is an $F$-dec-min element of $\odotZ{Q_{L}}$.  
\end{claim}

\Proof 
Suppose first that $z$ is an $F$-dec-min element of $\odotZ{Q}$.
Then $z$ is surely $F$-pre-dec-min in $\odotZ{Q}$ 
and hence, by Claim~\ref{Fdecmin}, $z$ is in $\odotZ{Q_{L}}$.  
If, indirectly, $\odotZ{Q_{L}}$ had an element $z'$ 
which is decreasingly smaller on $F$ than $z$, 
then $z$ could not have been an $F$-dec-min element of  $\odotZ{Q}$.

Conversely, let $z'$ be an $F$-dec-min element of $\odotZ{Q_{L}}$ and
suppose indirectly that $z'$ is not an $F$-dec-min element of $\odotZ{Q}$.  
Then any $F$-dec-min element $z$ of $\odotZ{Q}$ is 
decreasingly smaller on $F$ than $z'$.  
But any $F$-dec-min element of  $\odotZ{Q}$
is pre-dec-min on $F$ and hence, by Claim~\ref{Fdecmin}, 
$z$ is in $\odotZ{Q_{L}}$, contradicting the assumption that $z'$ 
is an $F$-dec-min element of $\odotZ{Q_{L}}$.  
\finbox

\medskip

For wider applicability, we must emphasize that, 
similarly to Theorem~\ref{MAIN}, the formulation of the next result relies on the
$0$-base-polyhedron $B$ itself, but not on the unique fully
supermodular function $p$ defining $B$.  
(The proof, however, does refer to $p$.)

\begin{theorem} 
\label{reduction} 
Let $D$, $(f,g)$, $F$, and $B$ be the same
as in Theorem~{\rm \ref{MAIN}}, and let $L:=L_\beta $ 
(defined in \eqref{(Lbetadef)}).  
Given \eqref{(hypo)}, there is a set
$F'\subset F$ for which an element $z$ of $\odotZ{Q}$ 
is an $F$-dec-min member of $\odotZ{Q}$ 
if and only if 
$z$ is an $F'$-dec-min member of \
$\odotZ{Q_{L}} =\odotZ {Q}(f_{L},g_{L};B_{L})$.  
In addition, the box $T(f_{L},g_{L})$ is narrow on $F-F'$ in the sense that 
$0\leq g_{L}(e)-f_{L}(e)\leq 1$ holds for every $e\in F-F'$.  
\end{theorem}

\Proof 
Apply Theorem~\ref{Lupmin} to $L=L_\beta $, and consider 
the face $B_{L}$ of $B$ ensured by the theorem.  
Let ${\cal C}_{L}$ be the chain describing $B_{L}$ 
(that is, 
$B_{L}=\{x\in B:  \widetilde x(C) = p(C)$ for each $C\in {\cal C}\}$ 
where $p$ denotes the fully supermodular function $p$ defining $B$).  
Let $(f_{L},g_{L})$ be the pair of bounding functions defined in \eqref{(f'g'def1)} 
and \eqref{(f'g'def2)}, and let 
$\odotZ{Q_{L}}:  = \odotZ{Q}(f_{L},g_{L};B_{L})$.

\begin{claim} 
\label{L'} 
The subset $L'\subseteq L$ consisting of the ${\cal C}_{L}$-entering elements 
of $L$ is non-empty.  
\end{claim}

\Proof 
Let $z$ be an element of $\odotZ{Q}$ which is pre-dec-min on $F$.  
By Claim~\ref{Fdecmin}, $z\in \odotZ{Q_{L}}$.  
By \eqref{(hypo)}, there is an edge $e$ in $F$, 
for which $z(e)=\beta =g(e)$, and hence $e\in L$.  
Since $g(e)=z(e)\leq g_{L}(e)\leq g(e)$ and $F$ contains no
$(f,g)$-tight edges, we have $f(e) < g(e) =g_{L}(e) =\beta $. 
This and definition \eqref{(f'g'def1)} imply that $e$ is ${\cal C}_{L}$-entering.
\finbox \medskip

Since $L'\not =\emptyset $ by the claim, we have that

\[ 
\hbox{ $F':=F-L'$ \ is a proper subset of $F$.}\ 
\]
\noindent 
We are going to show that $(f_{L},g_{L})$ and $F'$ meet 
the requirements of the theorem.  
Call two vectors in ${\bf Z}\sp A$ value-equivalent on $L'$ 
if their restrictions to $L'$ (that is, their projection to ${\bf Z}\sp {L'}$), 
when both arranged in a decreasing order, are equal.

\begin{claim} 
\label{valeq} 
The members of $\odotZ{Q_{L}}$ are value-equivalent on $L'$.  
\end{claim}

\Proof 
By Claim~\ref{Fdecmin}, the members of $\odotZ{Q_{L}}$ are exactly
those elements of $\odotZ{Q}$ which are pre-dec-min on $F$.  
Hence each member $z$ of $\odotZ{Q_{L}}$ has the same number $\mu$ 
of edges $e$ in $L$ for which $z(e)=\beta$.

As $F$ contains no $(f,g)$-tight edges, we have 
$z(e)\leq g_{L}(e)\leq \beta -1$ 
for every edge $e\in L-L'$ and hence each element $e$ of $L$
with $z(e)=\beta $ belongs to $L'$, 
from which
\[ 
 | \{e\in L':  z(e)= \beta \} | = \mu .
\]

Furthermore, we have $f_{L}(e)\geq \beta -1$ for every element $e$ 
of $L'$, from which $L'$ has exactly $ | L' | -\mu $ edges with
$z(e)=\beta -1$, implying that the members of $\odotZ{Q_{L}}$ are indeed
value-equivalent on $L'$.  
\finbox \medskip

Claim~\ref{f'g'} implies that the $F$-dec-min elements of $\odotZ{Q}$
are exactly the $F$-dec-min elements of $\odotZ{Q_{L}}$, and hence it
suffices to prove that an element $z$ of $\odotZ{Q_{L}}$ is an
$F$-dec-min member of $\odotZ{Q_{L}}$ 
if and only if 
$z$ is an $F'$-dec-min member of \ $\odotZ{Q_{L}}$.  
But this latter equivalence is an immediate consequence of Claim~\ref{valeq}.

To prove the last part of Theorem~\ref{reduction}, recall that
$F-F'=L'$ and $L'$ consisted of the ${\cal C}_{L}$-entering elements of $L$.  
But the definition of $(f_{L},g_{L})$ in \eqref{(f'g'def1)} implies
that $\beta -1\leq f_{L}(e)\leq g_{L}(e)=\beta $ for every element $e$ of
$L'$, that is, the box $T(f_{L},g_{L})$ is indeed narrow on $F-F'$.  
This completes the proof of Theorem~\ref{reduction}.
\BB

\subsection{Proof of the main result}

After these preparations, we are in a position to prove our main
result formulated in Section~\ref{Intro}.  
\medskip

\noindent 
{\bf Proof of Theorem~\ref{MAIN}}.  \ 
We use induction on $ | F | $. 
Since $f\sp{*}:=f$, $g\sp{*}:=g$, and $p\sp{*}:=p$
clearly meet the requirements of the theorem when $F=\emptyset $, 
we can assume that $F$ is non-empty.  
As before, we may assume that $F$ contains no $(f,g)$-tight edges.  
By Theorem~\ref{reduction}, 
it suffices to prove the theorem for $\odotZ{Q}(f_{L},g_{L};B_{L})$ and $F'$.
But this follows by induction since $F'$ is a proper subset of $F$.
\FBOX

\subsection{Graph orientations}

One of the starting points of the present investigation was the paper
of Borradaile et al.  \cite{BIMOWZ}, in which they considered 
(among others) strongly connected orientation of an undirected graph 
for which the indegree vector is decreasingly minimal 
(egalitarian in their term).  
They formulated a conjecture for characterizing dec-min strongly connected orientations.  
The conjecture has been proved in \cite{Frank-Murota.A} 
in a more general framework concerning $k$-edge-connected in-degree constrained orientations.  
The characterization immediately gave rise to an algorithm 
for computing a dec-min orientation in question.  
The approach of \cite{Frank-Murota.A}, however, does not say anything about strongly
connected dec-min orientations of a mixed graph.

It is known, however, (see, for example, \cite{Frank-book}) that even
the more general problem of finding a $k$-edge-connected in-degree
constrained orientation of a mixed graph can be formulated as a
special base-flow polyhedron problem (where the base-flow is defined
by a crossing supermodular function).  Therefore both Theorem
\ref{MAIN} and the algorithm developed for computing minimum cost
dec-min base-flow can be specialized to this dec-min orientation problem
on a mixed graph.

\section{Algorithm for minimizing the largest base-flow value on $F$}
\label{algo1}

Our next task is to describe a strongly polynomial algorithm to
compute the bounding pair $(f\sp{*},g\sp{*})$ and the face $B\sp{\triangledown}$
of $B$ in Theorem~\ref{MAIN}, the main result of the paper.  
In Section~\ref{decmin} we derived Theorem~\ref{MAIN} from Theorem~\ref{reduction}.  
This derivation showed that if an algorithm is available 
to compute the face $B_{L}$ of $B$ and the bounds $(f_{L},g_{L})$
occurring in Theorem~\ref{reduction}, then at most $| F |$
applications of this algorithm result 
in the requested $(f\sp{*},g\sp{*})$ 
and $B\sp{\triangledown}$.  
Therefore our main task is to show how the
bounds $(f_{L},g_{L})$ and face $B_{L}$ in Theorem~\ref{reduction} can be computed.

Theorem~\ref{reduction} itself was derived by applying Theorem~\ref{Lupmin} 
in the special case $L:=L_\beta $. 
Therefore, we can apply the algorithm outlined at the end of Section~\ref{minmaxL} 
once we are able to compute $\beta_F$ defined in \eqref{(betaF)}.  
Recall from Section~\ref{predecminF} that $\beta_F$ is nothing but 
the smallest integer for which $\odotZ{Q}$ has an element $z$ that is
$\beta_F$-covered on $F$, or in other words, every component of $z$
in $F$ is at most $\beta_F$.

As a preparation to computing $\beta :=\beta_F$, let $S$ be a finite ground-set, 
$h:2\sp{S}\rightarrow \underline{\bf Z}$ and 
$b:2\sp{S}\rightarrow \overline{\bf Z}_{+}$ 
be set-functions.  We call an integer $\mu$ \ {\bf good} 
with respect to $b$ and $h$ if $\mu b(X)\geq h(X)$ 
for every $X\subseteq S$.  
An integer that is not good is called {\bf bad}.  
The non-negativity of $b$ implies that if $\mu$ is good, 
then so is every integer larger than $\mu$. 
We assume that
\begin{equation} 
\mbox{$h(X)\leq 0$ \quad whenever \ $b(X)=0$},
\label{(vanjomu)} 
\end{equation} 
which is equivalent to requiring that there is a good $\mu$.  
We also assume that 
\begin{equation}
\mbox{there exists a subset \ $Y\subseteq S$ \ with \ $h(Y)>0$},
\label{(0isbad)} 
\end{equation} 
which is equivalent to requiring that the value $\mu =0$ is bad.  
Let $\mu_{\rm min}$ denote the smallest good integer.

\subsection{Computing the smallest good $\mu$} 
\label{SCpbmaxND}

We recall an algorithm from \cite{Frank-Murota.C} to compute $\mu_{\rm min}$.  
(In order to avoid confusion, here we use letter $h$ for
the function $p$ used in \cite{Frank-Murota.C}.)  
The number $\mu_{\rm min}$ is nothing but the maximum of 
\ $\lceil {h(X) / b(X)}\rceil $ \ 
over the subsets $X$ of $S$ with $b(X)>0$, and hence
the algorithm may be viewed as a variant of the Newton--Dinkelbach algorithm
\cite{GGJ17,Radzik}.

The algorithm works if a subroutine is available to 
\begin{equation}
\hbox{ find a subset $X\subseteq S$ 
maximizing $h(X) - \mu b(X)$ \ for any fixed integer $\mu \geq 0$.  }\ 
\label{(ND.routine)}
\end{equation} 
\noindent 
This routine will actually be needed only for special values of 
$\mu$ when $\mu =\lceil h(X)/\ell\rceil$ $\geq 0$
with $X\subseteq S$ and $1\leq \ell\leq M$, 
where $M$ denotes the largest finite value of $b$.

The algorithm starts with the bad $\mu_{0}:=0$.  
Let 
\[ 
X_{0} \in \arg\max \{ h(X)-\mu_{0}b(X) :  \ X\subseteq S \}, 
\] 
that is, \
$X_{0}$ is a set maximizing the function $h(X)-\mu_{0}b(X)=h(X)$.
Note that the badness of $\mu_{0}$ implies that $h(X_{0}) > 0$.
Since, by assumption, there is a good $\mu$, it follows that 
$\mu b(X_{0}) \geq h(X_{0})$
for some $\mu > 0$, and hence $b(X_{0})>0$.

The procedure determines one by one a series of pairs
$(\mu_{j},X_{j})$ for subscripts $j=1,2,\dots $, where each integer
$\mu_{j}$ is a tentative candidate for $\mu$ while $X_{j}$ 
is a non-empty subset of $S$ with $b(X_{j}) > 0$.  
Suppose that the pair
$(\mu_{j-1},X_{j-1})$ has already been determined for a subscript $j\geq 1$.  
Let $\mu_{j}$ be the smallest integer for which
$\mu_{j}b(X_{j-1})\geq h(X_{j-1})$, that is,
\[ 
\mu_{j}:  = \left\lceil 
\frac{h(X_{j-1})}{b(X_{j-1})} 
\right\rceil .
\]

If $\mu_{j}$ is bad, then let 
\[ 
X_{j} \in \arg\max \{ h(X)-\mu_{j}b(X) :\ X\subseteq S \}, 
\] 
that is, \ $X_{j}$ is a set
maximizing the function \ $h(X)-\mu_{j}b(X)$.  
Since $\mu_{j}$ is bad, we have $h(X_{j}) - \mu_{j}b(X_{j})>0$, 
which implies $b(X_{j}) > 0$ by the assumption \eqref{(vanjomu)}.

It was proved in \cite[Section 9]{Frank-Murota.C} (Claim 9.1) 
that if $\mu_{j}$ is bad for some subscript $j\geq 0$, then $\mu_{j} < \mu_{j+1}$.  
This implies that there is a first subscript $\ell\geq 1$ during the run of
the algorithm for which $\mu_{\ell}$ is good.  
Theorem 9.2 in \cite[Section 9]{Frank-Murota.C} 
states that $\mu_{\rm min}=\mu_{\ell}$, that is,
$\mu_{\ell}$ is the requested smallest good $\mu$-value and $\ell\leq M$.

\subsection{Computing $\beta_{F}$ in strongly polynomial time}
\label{SCminlargeflow}

In the previous subroutine to compute the smallest good integer,
$S$ denoted the ground-set of the occurring set-functions $h$ and $b$.
Now we turn to the problem of computing $\beta_F$ and apply the
subroutine to set-functions defined on ground-set $V$ in place of $S$.

As before, we suppose that there is an $(f,g)$-bounded base-flow 
(that is, $Q=Q(f,g;B)$ is non-empty), 
and also that $F$ contains no $(f,g)$-tight edges.  
Based on the algorithm in Section~\ref{SCpbmaxND}, 
we describe first a strongly polynomial algorithm to compute 
$\beta_{F}$ defined in \eqref{(betaF)}, 
which is the smallest integer for which $\odotZ{Q}$ has an element $z$ 
satisfying $z(e)\leq \beta_F$ for every edge $e\in F$.  
Note that $\beta_F$ can be interpreted as the smallest integer such that, 
by decreasing $g(e)$ to $\beta_F$ for each edge $e\in F$ with $g(e)>\beta_F$, 
the resulting $g\sp{-}$ and the unchanged $f$ continue 
to meet the inequality $f\leq g\sp{-}$ and 
the inequality $\varrho_{g\sp{-}}-\delta_{f}\geq p$ 
in \eqref{(baseflow-nonempty)}, which ensure that $Q(f,g\sp{-};B)$ is non-empty.

The first requirement for $\beta_F$ is that it should be 
at least the largest $f$-value on the edges in $F$, 
which is denoted by $f_{1}$.
Let $g_{1}>g_{2}>\cdots >g_{q}$ denote the distinct $g$-values of 
the edges in $F$, and let $L:=\{e\in F:  g(e)=g_{1}\}$.  
Let $\beta_{1}:=\max \{f_{1},g_{2}\}$.

It is known that there is a (purely combinatorial) strongly polynomial algorithm 
(see, for example \cite{FrankJ8,Frank-book}) 
to check whether a base-flow polyhedron is empty or not.  
We refer to such an algorithm as a base-flow feasibility subroutine.  
Such an algorithm either outputs an integral element of the base-flow polyhedron 
in question or else it outputs a set $X$ violating \eqref{(baseflow-nonempty)}.

With the help of this subroutine, we can check whether the $g$-value
$g_{1}$ on the elements of $L$ can uniformly be decreased to 
$\beta_{1}$ without destroying \eqref{(baseflow-nonempty)}.  
If this is the case, then either $\beta_{1}=f_{1}$, 
in which case a tight edge arises in $F$ and 
we can remove this tight edge from $F$, 
or $\beta_{1}=g_{2}$, 
in which case the number of distinct $g$-values becomes one smaller.  
Clearly, as the total number of distinct $g$-values in $F$ is 
at most $ | F | $, this kind of reduction may occur at most $ | F | $ times.

Therefore, we are at the case where $g_{1}$ cannot be decreased 
to $\beta_{1}$ without violating \eqref{(baseflow-nonempty)}, that is, 
$\beta_F > \beta_{1}$.  
We look for $\beta_F$ in the form 
$\beta_F=\beta_{1}+\mu $, 
that is, our goal is to compute the smallest positive integer $\mu$ 
such that $\odotZ{Q}$ has an element which is $(\beta_{1}+\mu )$-covered on $F$.  
To this end, we show that computing such a $\mu$ is nothing but 
finding a smallest good $\mu$ with respect to set-functions $h:=p'$ and $b$ 
to be defined as follows.

Recall that $L=\{e\in F:  g(e)=g_{1}\}$ and let $A_{0}:=A-L$ (that is,
$A_{0}$ is the complement of $L$ with respect to the whole edge-set $A$).  
Let $g'$ denote the function arising from $g$ by reducing
$g(e)$ on the elements of $L$ (where $g(e)=g_{1}$) to $\beta_{1}$.
Since $g'\geq f$ holds and hence $\varrho_{g'}-\delta_{f}$ is submodular, 
the set-function $p'$ on $V$ defined by 
\begin{equation}
\label{supmodfnforND} 
p'(Z):  = p(Z) - \varrho_{g'}(Z) + \delta_{f}(Z) 
\end{equation} 
is supermodular.  
Define a submodular function $b$ on $V$ by 
\begin{equation} \label{submodfnforND} b(Z):= \varrho_{L}(Z) . 
\end{equation} 
\noindent 
Recall that a non-negative integer $\mu$ is called good 
with respect to $b$ and $p'$ if 
\begin{equation}
\mu b(Z)\geq p'(Z) 
\label{(mubZ)} 
\end{equation} 
holds for every $Z\subseteq V$.

\begin{claim} 
\label{ekvimu} 
For a non-negative integer $\mu$, Condition \eqref{(baseflow-nonempty)} 
holds after increasing $g(e)=\beta_{1}$ uniformly by $\mu$ 
on the edges $e\in L$, or equivalently, 
\begin{equation} 
\mu \varrho_{L}(Z) + \varrho_{g'}(Z) - \delta_{f}(Z) \geq p(Z) \
\label{(mubZx)} 
\end{equation} 
holds for every $Z\subseteq V$ if and only if
$\mu$ is good with respect to $b$ and $p'$.  
\end{claim}

\Proof The claim follows immediately from the equivalence of
\eqref{(mubZ)} and \eqref{(mubZx)}.  
\FBOX

\medskip

Therefore our problem of computing $\beta_F$ is reduced to computing
the smallest good non-negative integer $\mu$. 
As $Q$ has no element which is $g_{1}$-covered on $F$, 
Claim~\ref{ekvimu} shows that $\mu =0$ is not good.  
Similarly, as $Q$ has an element which is $g_{1}$-covered on $F$, 
Claim~\ref{ekvimu} shows that $\mu :=g_{1}-\beta_{1}$ is good.

Since $b$ is submodular, $p'$ is supermodular, and we have 
$\max \{b(Z):Z\subseteq V\} \leq  | L | \leq  | A | $, 
we can apply the algorithm described in Section~\ref{SCpbmaxND} 
to $h:=p'$ and $b$.  
That algorithm needs the subroutine \eqref{(ND.routine)}
to compute a subset of 
$V$ maximizing $p'(Z) - \mu b(Z)$ \ ($Z\subseteq V$) 
for any fixed integer $\mu \geq 0$.  
This subroutine is applied at most $M$ times, 
where $M$ denotes the largest value of $b$.  
Since the largest value of $b$ is at most $ | A | $, 
the subroutine \eqref{(ND.routine)} is applied at most $ | A | $ times.
Furthermore, by the definition of $p'$ and $b$, the equivalent subroutine 
to minimize 
\[ 
\mu b(Z)-p'(Z)= \mu \varrho_{L}(Z) + \varrho_{g'}(Z) -\delta_{f}(Z) - p(Z) 
\] 
can be realized with the help of a submodular function minimizing algorithm.

Therefore, the smallest good $\mu$ can be computed in strongly
polynomial time, and hence the requested $\beta_F = \beta_{1}+ \mu $
may be assumed to be available.

\begin{remark} \rm \label{RMprincpat}
In the above argument, we have considered a supermodular function
$p'(Z) - \mu b(Z)$ with a parameter $\mu \geq 0$.
The principal partition is the central concept 
in a structural theory for 
submodular/supermodular functions
that considers a parametrized family of 
supermodular functions
of the form 
$\hat p(Z) - \alpha \hat b(Z)$,
where 
$\hat p$ is a submodular function,
$\hat b$ is a submodular function,
and $\alpha$ is a nonnegative real-valued parameter.
See 
Iri \cite{Iri79} for an early survey and 
Fujishige \cite{Fuj09bonn} for a comprehensive historical and technical account.
As was pointed out by Fujishige \cite{Fujishige80},
the dec-min element of a base-polyhedron in continuous variables
is closely related to the principal partition 
in the special case when 
$\hat b(Z) = |Z|$.
In addition, it is shown recently \cite{FMdecminRZ} that 
the discrete dec-min problem on 
an integral base-polyhedron
(that is, on an M-convex set)
is also closely related to the principal partition
where the parameter $\alpha$ is restricted to integers.
It would be an interesting research topic
to investigate the relationship between the dec-min 
base-flow problem (continuous and discrete)
with the principal partition.
In this connection, it is noted that 
the special case of the dec-min mod-flow problem
where the specified set $F$ of edges is the star leaving a single (source) node
is a special case of the dec-min problem on a base-polyhedron (continuous and discrete).
$\bullet $ 
\end{remark}

\section{Existence of an $F$-dec-min base-flow} \label{vandecmin}

In the previous sections, we assumed that the bounding functions $f$
and $g$ are finite-valued on $F$ because this certainly ensured the
existence of an $F$-dec-min element of $\odotZ{Q}$ where $Q=(f,g;B)$ is
a non-empty base-flow polyhedron.  
The goal of the present section is to describe a characterization 
of the existence of $F$-dec-min elements when there are 
no a priori assumptions on the finiteness of $f$ and $g$.  
Since every member of $\odotZ{Q}$ is trivially $F$-dec-min when $F=\emptyset $, 
we assume that the specified subset $F$ of $A$ is non-empty.

First, we exhibit an easy reduction by which $g$ can be made finite-valued on $F$ 
without changing the set of $F$-dec-min elements of $\odotZ{Q}$.

\begin{lemma} 
\label{gfinonF} 
There is a function $g\sp{-} \ (\leq g)$ on $A$
which is finite-valued on $F$ such that the (possibly empty) set of
$F$-dec-min elements of $\odotZ{Q}:=\odotZ{Q}(f,g;B)$ 
is equal to the set of $F$-dec-min elements of 
$\odotZ{Q'}:=\odotZ{Q}(f,g\sp{-};B)$.
\end{lemma}

\Proof 
Let $z_{1}$ be an element of $\odotZ{Q}$ and let $\beta $ denote
the maximum value of its components.  
Define $g\sp{-}$ as follows:
\begin{equation} 
g\sp{-}(e):= \begin{cases} 
\min \{g(e),\beta \} & \quad \hbox{if}\ \  \ e\in F, \ \cr 
g(e) & \quad \hbox{if}\ \ \ e \in A-F.  
\end{cases}
\label{(g'def)} 
\end{equation}

As $g\sp{-}\leq g$, we have $\odotZ{Q'} \subseteq \odotZ{Q}$.
In particular, an $F$-dec-min element $z'$ of $\odotZ{Q'}$ is in $\odotZ{Q}$, 
and we claim that $z'$ is actually $F$-dec-min in $\odotZ{Q}$, as well.  
Indeed, if we had an element $z'' \in \odotZ{Q}$ which is
decreasingly smaller on $F$ than $z'$, then $z''$ is not in $\odotZ{Q'}$, 
that is, $z''$ is not $(f,g\sp{-})$-bounded.  
Therefore there is an edge $a\in F$ for which $z''(a)>\beta $, implying that 
$\max \{z''(e) :  e\in F\} >\beta \geq \max \{z'(e):e\in F\}$.  
But this contradicts the assumption that $z''$ is decreasingly smaller on $F$ than $z'$.

Conversely, suppose that $z$ is an $F$-dec-min element of $\odotZ{Q}$.  
Since the largest component of $z_{1}$ is $\beta $, 
the largest component of $z$ is at most $\beta $, and hence $z\in \odotZ{Q'}$.
This and $\odotZ{Q'} \subseteq \odotZ{Q}$ imply that $z$ is an
$F$-dec-min element of $\odotZ{Q'}$.  
\finbox \medskip

Let 
\begin{equation} 
{\cal P}:=\{Z\subseteq V:  \ \min \{\widetilde x(Z):  x\in B\} >-\infty \}.  
\label{(pveges)} 
\end{equation} 
\noindent 
Note that if $p$ denotes the unique fully supermodular function defining $B$, 
then 
\[ 
{\cal P}=\{Z\subseteq V:  \ p(Z)>-\infty \}.  
\] 
This implies that ${\cal P}$ is closed under taking intersection and union. 
 We also have $V\in {\cal P}$ since $B$ is a $0$-base-polyhedron.  
Therefore, each element $u\in V$ is contained 
in the unique smallest member $P(u)$ of ${\cal P}$
(which is the intersection of all members of ${\cal P}$ containing $u$).

To formulate the main result of this section, we introduce a set $J$
of edges that encodes ${\cal P}$.  
Let
\begin{equation} 
J:=\{uv:  \ v\in P(u)-u\}.  
\label{(Jdef)} 
\end{equation}
\noindent 
We refer to the elements of $J$ as {\bf jumping edge}s.
Clearly, $e=uv$ is a jumping edge 
precisely if $e$ does not leave any member of ${\cal P}$, 
and a subset $X\subseteq V$ is a member of ${\cal P}$ 
precisely if $\delta_J(X)=0$.

Let 
\[ 
A_{1}:=\{e:  e\in A, \ f(e)=-\infty \},
\qquad
A_{2}:=\{vu:  uv\in A-F, \ g(uv)=+\infty \},
\] 
and let 
\begin{equation} 
A\sp{*} := J \ \cup \ A_{1} \ \cup \ A_{2}.  
\label{(Ainfdef)} 
\end{equation}

\begin{theorem} 
\label{findecmin} 
Let $D=(V,A)$ be a digraph and $F\subseteq A$ a non-empty subset of edges, 
and let $B$ be an integral $0$-base-polyhedron.  
Let $f:A\rightarrow \underline{\bf Z}$ and
$g:A\rightarrow \overline{\bf Z}$ 
be bounding functions on $A$ with $f\leq g$ 
such that the base-flow polyhedron $\odotZ{Q}:=\odotZ{Q}(f,g;B)$ is non-empty.  
Then there exists an $F$-dec-min element of $\odotZ{Q}$
if and only if 
there is no di-circuit $C$ with $C\cap F\not =\emptyset$ 
in the digraph $D\sp {*}=(V,A\sp{*})$ defined by \eqref{(Ainfdef)}.
\end{theorem}

\Proof 
Let $p:2\sp V\rightarrow \underline{\bf Z}$ be 
the unique fully supermodular function on $V$ determining $B$, 
that is, $B=B'(p)$.
Since $B$ is a $0$-base-polyhedron, $p(V)=0$.  
Suppose first that $D\sp{*}$ includes a di-circuit $C$ intersecting $F$, 
and assume, indirectly, that there exists an $F$-dec-min member $z$ of $\odotZ{Q}$.
For $uv\in A$, define $z'(uv)$ as follows:
\begin{equation} 
z'(uv):= \begin{cases} 
z(uv)-1 & \quad \hbox{if}\ \ \ uv\in C\cap A_{1}, \cr 
z(uv)+1 & \quad \hbox{if}\ \ \ vu \in C\cap A_{2}, \cr 
z(uv) & \quad\hbox{otherwise}.  \end{cases} 
\label{(z'def)} 
\end{equation}

Let $Z$ be any subset of $V$.  
Clearly, 
\begin{equation} 
\varrho_C(Z) = \delta_C(Z).  
\end{equation} 
Let $C_{1}:=C\cap A_{1}$, $C_{2}:=C\cap A_{2}$, and $C_J:=C\cap J$.  
Then 
\begin{equation} 
  \varrho_{C_{1}}(Z) + \varrho_{C_{2}}(Z) + \varrho_{C_J}(Z) 
  = \delta_{C_{1}}(Z) + \delta_{C_{2}}(Z) + \delta_{C_J}(Z)
\label{(rodelta1)} 
\end{equation}
and hence 
\begin{equation} 
 \varrho_{C_J}(Z) - \delta_{C_J}(Z) 
 = \delta_{C_{1}}(Z) + \delta_{C_{2}}(Z) - \varrho_{C_{1}}(Z) - \varrho_{C_{2}}(Z).  
\label{(rodelta2)} 
\end{equation} 
Furthermore, 
\begin{align} 
 \varrho_{z'}(Z) &= \varrho_{z}(Z) - \varrho_{C_{1}}(Z) + \delta_{C_{2}}(Z),
\label{(ro)}
\\
\delta_{z'}(Z) &= \delta_{z}(Z) - \delta_{C_{1}}(Z) + \varrho_{C_{2}}(Z), 
\label{(delta)} 
\end{align} 
\noindent 
 from which, by recalling the notation 
$\Psi_{z}:=\varrho_{z}-\delta_{z}$, 
we get 
\begin{equation} 
\Psi_{z'}(Z) 
= \Psi_{z}(Z) - \varrho_{C_{1}}(Z) + \delta_{C_{2}}(Z) + \delta_{C_{1}}(Z) - \varrho_{C_{2}}(Z) 
= \Psi_{z}(Z) + \varrho_{C_J}(Z) - \delta_{C_J}(Z).  
\label{(psiz')} 
\end{equation}

\begin{claim} 
\label{z'also} 
The vector $z'$ defined in \eqref{(z'def)} is also in $\odotZ{Q}$.  
\end{claim}

\Proof 
It follows immediately from the definitions of $D\sp{*}$ and $z'$ 
that $f\leq z'\leq g$.  
To prove that $\Psi_{z'}=\varrho_{z'}-\delta_{z'}\geq p$, 
let $Z$ be a member of ${\cal P}$.  
By the definition of $J$, no jumping edges leave $Z$, and hence 
$\delta_{C_J}(Z)=0$.  From \eqref{(psiz')} we get 
\[
 \Psi_{z'}(Z) = \Psi_{z}(Z) + \varrho_{C_J}(Z) \geq p(Z) + 0 = p(Z),
\] 
as required.  
\FBOX

\medskip

By Claim~\ref{z'also}, $z'\in \odotZ{Q}$.  
Since $F\cap C\not =\emptyset $, $z'$ is decreasingly smaller on $F$ than $z$,
 a contradiction, 
showing that in the present case no dec-min element of $\odotZ{Q}$ can exist.

\medskip

To see the converse, suppose that there is no di-circuit of $D\sp{*}$ intersecting $F$.  
We want to prove that there is an $F$-dec-min element of $\odotZ{Q}$.

We claim that it suffices to prove this statement in the special case
when $g$ is finite-valued on $F$. 
Indeed, consider the function $g\sp{-}$ introduced in \eqref{(g'def)}.  
As $g\sp{-}\leq g$, there is no di-circuit described in the theorem 
with respect to $(f,g\sp{-})$.  
Now $g\sp{-}$ is finite-valued on $F$, and if there exists 
an $F$-dec-min $(f,g\sp{-})$-bounded base-flow $z$, then it follows from 
Lemma~\ref{gfinonF} that $z$ is an $F$-dec-min $(f,g)$-bounded base-flow.

Therefore, we can assume that $g$ is finite-valued on $F$.  
In this case, 
\[
A_{2}= \{vu:  uv\in A, g(uv)=+\infty \}.
\]

\begin{claim} 
\label{l-bound.e0} 
Let $S\subset V$ be a set for which 
$\delta_{A\sp{*}}(S)=0$, 
and let $e_{0}\in F$ be an edge entering $S$.  
Then, for any $z\in \odotZ{Q}$, we have
\begin{equation} 
z(e_{0}) \geq p(S) - [ \varrho_{g}(S) - g(e_{0}) ] + \delta_{f}(S),
\label{(zlbound)} 
\end{equation} 
and the right-hand side is finite.  
\end{claim}

\Proof 
Since $z\leq g$ and $e_{0}$ enters $S$, we have 
\[
\varrho_{z}(S) - z(e_{0}) \leq \varrho_{g}(S) - g(e_{0}),
\] 
 from which 
\[
p(S) \leq \varrho_{z}(S)- \delta_{z}(S)
= z(e_{0}) + [\varrho_{z}(S) - z(e_{0})] - \delta_{z}(S) 
\leq z(e_{0}) + [ \varrho_{g}(S) - g(e_{0})] - \delta_{f}(S),
\]
implying \eqref{(zlbound)}.  
Furthermore, $\delta_{A\sp{*}}(S)=0$
implies that $f(e)>-\infty $ for every edge $e$ of $D$ leaving $S$ and
that $g(e) <+\infty $ for every edge $e$ of $D$ entering $S$.  
$\delta_{A\sp{*}}(S)=0$ also implies that no jumping edge leaves $S$, which is
equivalent to saying that $P(u)\subseteq S$ for each $u\in S$.  
But this latter property implies that $S\in {\cal P}$, that is, $p(S)$ is
finite, from which the finiteness of the right-hand side of
\eqref{(zlbound)} follows.  
\finbox \medskip

Assume indirectly that $\odotZ{Q}$ has no $F$-dec-min element, that is,
for every element of $\odotZ{Q}$ there exists another one which is
decreasingly smaller on $F$.  
This implies that there is an edge
$e_{0}=ts$ in $F$ 
such that, for an arbitrary small integer $K$, there is 
an element $z$ of $\odotZ{Q}$ satisfying $z(e_{0})\leq K$.

\begin{claim} 
\label{vanut} 
There exists an $st$-dipath $P$ in $D\sp{*}$.
\end{claim}

\Proof 
Suppose, indirectly, that the set $S$ of nodes reachable from $s$ 
in $D\sp{*}$ does not contain $t$.  
Since no edge of $D\sp{*}$ leaves $S$ and $e_{0}$ enters $S$, 
it follows from Claim~\ref{l-bound.e0} 
that there is a finite lower bound for $z(e_{0})$, a contradiction.  
\finbox \medskip

The di-circuit formed by $e_{0}=ts$ and the $st$-dipath $P$ ensured by
Claim~\ref{vanut} meets the requirement of the theorem.  
\BB

\medskip \medskip

\noindent {\bf Extension of Theorem~\ref{MAIN}} \ 
With the help of
Theorem~\ref{findecmin} and Lemma~\ref{gfinonF}, Theorem~\ref{MAIN}
can be extended to the case when $(f,g)$ is not assumed to be
finite-valued on $F$, but only the existence of 
a di-circuit in $D\sp{*}$ intersecting $F$ is excluded 
(which is equivalent, by Theorem~\ref{findecmin}, 
to the existence of an $F$-dec-min member of $\odotZ{Q}$).

\begin{theorem}
\label{MAINb} 
Let $D=(V,A)$ be a digraph endowed with
integer-valued lower and upper bound functions 
$f:A\rightarrow \underline{\bf Z}$ 
and $g:A\rightarrow \overline{\bf Z}$ 
for which $f\leq g$.  
Let $B$ be an integral $0$-base-polyhedron 
for which the base-flow polyhedron $Q=Q(f,g;B)$ is non-empty.  
Let $F\subseteq A$ be a specified subset of edges 
for which there exists an $F$-dec-min element of $\odotZ{Q}$.  
Then there exists a face $B\sp{\triangledown}$ of $B$
and there exists a pair $(f\sp{*},g\sp{*})$ 
of integer-valued bounding functions on $A$ 
with $f\leq f\sp{*} \leq g\sp{*} \leq g$ 
such that an element $z\in \odotZ{Q}$ is $F$-dec-min 
if and only if 
$z\in \odotZ{Q}(f\sp{*},g\sp{*};B\sp{\triangledown})$.  
Moreover, $0\leq g\sp{*}(e)-f\sp{*}(e)\leq 1$ for every $e\in F$.  
\end{theorem}

\Proof 
By Lemma~\ref{gfinonF}, we can assume that $g$ is finite-valued on $F$.  
By Theorem~\ref{findecmin}, there is no di-circuit $C$ in
$D\sp{*}$ with $C\cap F\not =\emptyset $, implying that, for every edge
$e=ts\in F$, the set $S_e$ reachable in $D\sp{*}$ from $s$ meets 
the inequality \eqref{(zlbound)} for any member $z$ of $\odotZ{Q}$.  
As the right-hand side of \eqref{(zlbound)} is finite by Claim~\ref{l-bound.e0}, 
there is a finite lower bound 
\begin{equation} 
f'(e):= p(S_e) - [\varrho_{g}(S_e) - g(e) ] + \delta_{f}(S_e) 
\label{(f'(e))} 
\end{equation} 
for $z(e)$.  
In this way, each $(-\infty )$-valued lower bound on the edges in $F$ 
can be made finite, and the original Theorem~\ref{MAIN} applies.  
\finbox

\medskip

We emphasize that for each $e\in F$ the set $S_e$ occurring in the proof 
is easily computable and hence so is the finite lower bound
$f'(e)$ given in \eqref{(f'(e))}.  
Therefore this reduction to the case
when $(f,g)$ is finite-valued on $F$ is algorithmic.

\begin{remark} \rm
Since the intersection $Q$ of two g-polymatroids is known 
\cite{FrankP6}
to be a base-flow polyhedron, 
Theorem~\ref{findecmin} can be specialized to this case. 
Furthermore, Theorem~\ref{MAINb} implies a slight extension of Corollary~\ref{M2}
when, instead of assuming the boundedness of $Q$, 
we assume only the existence of an $F$-dec-min element of $\odotZ{Q}$.
$\bullet $ 
\end{remark}

\section{Fractional dec-min elements}

Throughout the paper we concentrated exclusively 
on integral dec-min elements, but analogous questions 
concerning fractional dec-min elements of a base-flow polyhedron also make sense.  
In \cite{Frank-Murota.C}, we proved the following proposition.

\begin{proposition} 
If a convex subset $P$ of ${\bf R} \sp n$ admits a dec-min element $x$, 
then $x$ is the unique dec-min element.  
\end{proposition}

This theorem can be applied to base-flow polyhedra, and 
since base-flow polyhedra are closed under projection, we even have 
the following.

\begin{claim} 
Let $Q=Q(f,g;B)$ be a non-empty base-flow polyhedron defined on
the edge-set of digraph $D=(V,A)$, and let $F\subseteq A$ be a subset of edges.  
If $Q$ admits a (fractional) element $x$ which is $F$-dec-min, 
then the restriction $x' | F$ of every $F$-dec-min element $x'$ of $Q$ to $F$ 
is the same as the restriction $x | F$ of $x$ to $F$.  
\end{claim}

Concerning the existence of an $F$-dec-min element of $Q$, we have the
following theorem, which is the continuous counterpart of Theorem~\ref{findecmin}.

\begin{theorem} 
\label{findecmin-R} 
There exists an $F$-dec-min element of $Q$ 
if and only if 
there is no di-circuit $C$ with $C\cap F\not =\emptyset $ 
in the digraph $D\sp {*}=(V,A\sp{*})$ defined by \eqref{(Ainfdef)}.  
\end{theorem}

\Proof 
The proof is essentially the same as that of Theorem~\ref{findecmin}, 
except that the definition of $z'(uv)=z(uv) \pm 1$ in \eqref{(z'def)} 
should be changed to $z'(uv) := z(uv) \pm \delta$
using a sufficiently small $\delta > 0$ to meet the condition $z' \in Q$.  
\finbox \medskip

It remains to be a task for future research to construct an algorithm for
computing an $F$-dec-min element of $Q$, when it exists.



\paragraph{Acknowledgement} \ 
The research was partially supported by the 
Hungarian Scientific Research Fund - OTKA, No. NKFIH-FK128673
and by JSPS KAKENHI Grant Number JP20K11697. 
The authors appreciate insightful comments of a referee
that led to Remarks \ref{RMconvmin} and \ref{RMprincpat}.








\newpage

\tableofcontents 

\newpage


\end{document}